\documentclass[11pt, a4paper]{article}

\usepackage[applemac]{inputenc}
\usepackage{color}
\usepackage{amsmath,amssymb,amsthm,amscd}

\usepackage{amsmath,amsfonts,amsthm}
\newtheorem{thm}{Theorem} [section]
\newtheorem{cor}[thm]{Corollary}
\newtheorem{lem}[thm]{Lemma}
\newtheorem{prop}[thm]{Proposition}
\theoremstyle{definition}
\newtheorem{defn}[thm]{Definition}

\newtheorem{claim}{Claim}
\newcommand{\G}{\Gamma}
\renewcommand{\>}{\rangle}

\newcommand{\C}{\mathcal C}
\newcommand{\ap}{\textquoteleft}
\newcommand{\ch}{\textquoteright}
\newcommand{\q}[1]{\ensuremath{\langle#1\rangle}}
\newcommand{\sq}[2]{\ensuremath{\langle#1\vert#2\rangle}}
\newcommand{\lh}[1]{\ensuremath{\left|#1\right|}}
\newcommand{\restr}{\upharpoonright}
\renewcommand{\epsilon}{\varepsilon}
\date{}
\begin{document}
\author{A. Ould Houcine,  D.  Vallino}
\title{Algebraic \& definable closure in free groups}
\maketitle
\begin{abstract} We study  algebraic closure and its relation with definable closure in free groups and  more generally in torsion-free hyperbolic groups. Given a torsion-free hyperbolic group $\G$ and a  nonabelian subgroup $A$ of $\G$, we describe $\G$ as a constructible  group from the algebraic closure of $A$ along cyclic subgroups.  In particular, it follows  that the algebraic closure of $A$ is finitely generated, quasiconvex and hyperbolic.

Suppose that $\G$ is free. Then  the definable closure of $A$ is a free factor of the algebraic closure of $A$ and the rank of these groups is bounded by that of $\G$.  We prove  that the algebraic closure of $A$ coincides with  the vertex group containing $A$ in the generalized malnormal cyclic JSJ-decomposition of $\G$ relative to $A$.  If the rank of  $\G$ is bigger than  $4$,  then $\G$ has  a subgroup $A$ such that the definable closure of $A$ is a proper subgroup of the algebraic closure of $A$. This answers a question of  Sela.  \end{abstract}


\section{Introduction}

In field theory, an element $b$ is called algebraic over a field $K$ if it is a root of some non-zero polynomial with coefficients in $K$. This notion is very  fruitful and has many applications in mathematics, as Galois theory shows. Its analogues in more general contexts were extensively  studied. Model theory generalizes the notion as follows. Given a model $\mathcal M$, in a first-order  language $\mathcal L$,  and a subset $A$ of $\mathcal M$, an  element $b$ is said to be  \textit{algebraic} over $A$, if there exists an $\mathcal L$-formula $\varphi(x)$, with parameters from $A$, such that $\mathcal M$ satisfies $\varphi(b)$ and  the set $\{c \in \mathcal M| \mathcal M \models \varphi(c)\}$ is finite.  The \textit{algebraic closure} of $A$, denoted $acl(A)$,  is the set of algebraic elements over $A$.   If $\{c \in \mathcal M| \mathcal M \models \varphi(c)\}$ is a singleton, then $b$ is said to be \textit{definable} over $A$, and one defines analogously the \textit{definable closure} of $A$, denoted  $dcl(A)$, as the set of definable elements over $A$.

It is well-known,  in the context of algebraically closed fields,  that the above model-theoretic  notion coincides with the usual one by using the quantifier elimination theorem of Tarski; i.e. $b$ is algebraic over $K$ (in the sense of the theory of fields) if and only if $b \in acl(A)$ (see for instance \cite[Proposition 3.2.15]{Marker}).

Algebraic closure plays an important role in the study of strongly minimal theories and more generally finite dimensional and stable theories. For instance it permits to define, in a suitable context, Zariski's geometries.  It is also an essential piece in the study of model-theoretic Galois theory.  Poizat has developed a Galois theory for theories which eliminate  imaginaries \cite{poizat-imaginaire}, and  Casanovas and Farr\'e studied  degree of elimination of imaginaries needed to have a Galois correspondance \cite{Casanovas-weak}. More recently, Medvedev and Takloo-Bighash have carried out some notions of  Galois theory in the setting of first-order theories \cite{medvedev-takloo}.

Sela  has shown that free groups and more generally torsion-free hyperbolic groups are stable \cite{Sela-stab}. He has also shown a geometric  elimination of imaginaries in torsion-free hyperbolic groups \cite{sela-imaginaries}. This can be certainly used to develop Galois theory of free groups. Miasnikov, Ventura and Weil have developed algebraic extensions in free groups \cite{miasnikov-ventura-weil-algebraic}, which correspond essentially to the notion of algebraic closure defined above but restricted to quantifier-free formulas.

In 2008, Sela asked, given a free group $F$ of finite rank and a subset $A$ of $F$, if  the algebraic and the definable closure of $A$ coincide.  In this paper  we study the algebraic and the definable closure in free groups. In particular we give  a negative answer  to the question of Sela for  free groups of rank $\geq 4$  and a positive answer for   free groups of  rank 2.

It is rather easy to see that $acl(A)$ and $dcl(A)$ are $\mathcal L$-substructures of $\mathcal M$, and in particular,  when $\mathcal M$ is a group,   they are subgroups. As usual, to axiomatize  group theory,  we use the language $\mathcal L=\{., ^{-1}, 1\}$, where $.$ is interpreted as multiplication, $^{-1}$ is interpreted as the function which sends every element to its inverse and $1$ is interpreted as the trivial element. Let $\Gamma$ be a group and $A$ a subset of $\Gamma$. It is not hard to see that $A$ and the subgroup generated by $A$ have the same algebraic closure; similarly for the definable closure. Hence, without loss of generality we may assume that $A$ is a subgroup.  We note also that  if $\G$ is torsion-free and hyperbolic and if $A$ is nontrivial and abelian, then the algebraic closure  and the definable closure of $A$ coincide with the centralizer of $A$ (see Lemma \ref{acl-abelian}).

The main results of this paper are as follows.  One of the first natural questions is to see the constructibility  of $\G$ from the algebraic closure.

\begin{thm}\label{thm-princip-1} Let $\G$ be a torsion-free hyperbolic group and $A$ a nonabelian subgroup of $\G$. Then $\G$ can be constructed from $acl(A)$ by a finite sequence of amalgamated free products and HNN-extensions along cyclic subgroups. In particular,  $acl(A)$ is finitely generated, quasiconvex and hyperbolic.
\end{thm}

In geometric group theory, given a finitely generated group $\G$ and a set $\mathcal C$ of subgroups of $\G$, one studies  the link between the various possible graph of groups decompositions   of $\G$,  with edge groups from $\mathcal C$ (i.e. splittings of $\G$ over $\mathcal C$). Grushko and  Kurosh showed  that there is a canonical \textit{free} decomposition (i.e. with trivial edge groups)  from which all other free decompositions can be obtained by some particular operations. At this point, it becomes natural to seek similar canonical splittings for larger classes of groups $\mathcal C$.

Roughly speaking a JSJ-decomposition of $\G$ over $\mathcal C$ is  a \textit{canonical} graph of groups decomposition of $\G$ over $\mathcal C$, from which all other splittings of $\G$ over $\mathcal C$  can be obtained through some natural  operations. The uniqueness of such a decomposition is not generally guaranteed, but all these decompositions share the   most important necessary properties.

The theory of JSJ-decompositions has its origin in the work of Johannson, and Jaco and Shalen, who developed a theory of cutting irreducible three-dimensional manifolds into pieces along tori and annuli \cite{Jaco-Shalen}.  One can describe such decompositions in terms of splittings of the relevant fundamental group. A group theoretic version was developed by Kropholler  \cite{kropholler}. Later Sela constructed JSJ-decompositions for torsion-free hyperbolic groups over cyclic subgroups \cite{Sela-hopf} and then Sela and Rips \cite{sela-JSJ} extended it to general torsion-free finitely presented groups. Other constructions of JSJ-decompositions  for various classes of groups $\mathcal C$ have been carried out by many authors.

JSJ-decompositions have  many applications and were successfully used by Sela to solve the isomorphism problem of torsion-free hyperbolic groups and to develop  diophantine geometry over  free (and hyperbolic)  groups  in the solution of   Tarski's conjecture.

The following theorem connects the notions of  algebraic closure and cyclic  JSJ-decompositions in  free groups. For the precise notions of  JSJ-decompositions which we  use, we  refer  the reader at the end of Section 2.

\begin{thm}\label{thm-princip-3} Let $\G$ be a free group of finite rank and let $A$ be a nonabelian subgroup of $\G$.  Then $acl(A)$  coincides with  the vertex group containing $A$ in the generalized malnormal cyclic JSJ-decomposition of $\G$ relative to $A$.
\end{thm}

Strictly speaking the notion of \textit{JSJ-decompositions} used in the previous theorem is not a JSJ-decomposition in the sense of \cite{Guirardel-Levitt-JSJ1}. However it possesses the most important properties of  JSJ-decompositions of \cite{Guirardel-Levitt-JSJ1}. By using the definition given in \cite{Guirardel-Levitt-JSJ1}, the conclusion of the previous theorem is the following : $acl(A)$ coincides with the \textit{elliptic abelian  neighborhood} of the vertex group containing $A$ in  the cyclic JSJ-decomposition of $\G$ relative to $A$, where we suppose that $\G$ is freely indecomposable relative to $A$.

   We will also be interested in a restricted notion of the algebraic closure. Given a group $\G$ and a subgroup $A$,  the \textit{restricted algebraic closure}, denoted $racl(A)$,  is defined as follows. An  element $\gamma$ is in  $racl(A)$ if and only if   its orbit $\{f(\gamma) | f \in Aut(F/A)\}$ is finite, where $Aut(F/A)$ is  the group of  automorphisms of $F$ fixing $A$ pointwise.  Note that $racl(A)$ is a subgroup and contains $acl(A)$.   It turns out that, when $\G$ is a torsion-free hyperbolic group and $A$ is nonabelian,  $racl(A)$ coincides with  the vertex group containing $A$ in the generalized malnormal cyclic JSJ-decomposition of $\G$ relative to $A$ (see Proposition \ref{thm-place-of-acl-in-splitting}). Similarly here by using the definition of JSJ-decompositions  of  \cite{Guirardel-Levitt-JSJ1}, the conclusion is that $racl(A)$ coincides with   the \textit{elliptic abelian  neighborhood} of the vertex group containing $A$ in the  cyclic JSJ-decomposition of $\G$ relative to $A$.  Theorem \ref{thm-princip-3} shows that in free groups, we get an  identity between  restricted  and  algebraic closure.

Notice that, as a corollary of the general version of Theorem \ref{thm-princip-1} (see Corollary \ref{cor-const}), we  have the following.  If  $\G$ is  a free group of finite rank and $A$ is a nonabelian subgroup of $\G$, then the rank of $acl(A)$ is bounded by the rank of $\G$. In fact, we will show that if $acl(A) \leq K\leq \G$, where $K$ is finitely generated, then $rk(acl(A))\leq rk(K)$; that is $acl(A)$ is compressed in the sense of \cite{ventura-martino-compress}.

Regarding the relation between  algebraic  and  definable closure, though generally  they are   different, at least we can  assert the following.

\begin{thm}\label{thm-princip-4} Let $\G$ be a free group of finite rank and $A$ a nonabelian subgroup of $\G$.  Then $dcl(A)$ is a free factor of $acl(A)$.
\end{thm}

Combining this with Lemma \ref{acl-abelian} below, it follows that when the rank of $\G$ is two, then $acl(A)=dcl(A)$ for any nontrivial subgroup $A$ of $\G$. However, this is not true in higher rank free groups.

\begin{thm} \label{thm-princip-5} Any free group $\G$ of rank  $n \geq 4$ can be written as an HNN-extension $\G=\<H, t| u^t=v\>$,   such that $H$ has a proper  subgroup $A$ with $acl(A)=H$ and $dcl(A)=A$.
\end{thm}

This paper is organized as follows. In  next section we recall the material that
we require around notions in model theory,  $\G$-limit groups  and the tools needed in the sequel.   Section 3 concerns  constructibility and its main purpose is the proof of Theorem \ref{thm-princip-1}.  The proof of that theorem follows the same strategy as the one used by Sela to prove constructibility of limit groups; however we need to analyze the place of  algebraic closure more carefully. Section 4 is devoted to the study of the place of  algebraic closure in the JSJ-decomposition and  we show Theorem  \ref{thm-princip-3}.  Section 5  deals with the proofs of Theorems \ref{thm-princip-4} and \ref{thm-princip-5}.

\bigskip

\noindent\textbf{Acknowledgements}.  The first named  author wishes to express his thanks to Z. Sela  for suggesting the problem.

\section{Prerequisites}

The aim of this section is to give the background needed in the sequel. The first subsection deals with notions from model theory; for more details the reader is referred to \cite{Hodges(book)93, Marker}. Notions around limit groups and abelian JSJ-decompositions are exposed in the second subsection.

\subsection{Model theory}

Given a language $\mathcal L$, an $\mathcal L$-structure $\mathcal M$  and an $\mathcal L$-formula $\varphi (\bar x)$, where $\bar x$ is a tuple of variables of length $n$,  we denote by $\varphi(\mathcal M)$ the set $\{\bar m \in \mathcal M^n | \mathcal M \models \varphi(\bar m)\}$.  Let $\mathcal M$ be an $\mathcal L$-structure   and $A$ a subset of $\mathcal M$. The \textit{algebraic closure} (resp. \textit{existential algebraic closure})  of $A$, denoted  $acl_{\mathcal M}(A)$ (resp. $acl_{\mathcal M}^{\exists}(A)$),  is the set of elements $x \in \mathcal M$  such that there exists a $\mathcal L$-formula  (resp. an existential $\mathcal L$-formula)  $\phi(x)$ with parameters from $A$  such that $\mathcal M \models \phi(x)$ and $\phi(\mathcal M)$ is finite.  The \textit{definable closure} (resp. \textit{existential definable closure})  of $A$, denoted  $dcl_{\mathcal M}(A)$ (resp. $dcl_{\mathcal M}^{\exists}(A)$),  is the set of elements $x \in \mathcal M$  such that there exists a formula  (resp. an existential formula)  $\phi(x)$ with parameters from $A$  such that $\mathcal M \models \phi(x)$ and $\phi(\mathcal M)$ is a singleton.   The previous notions are connected to  other notions of closedness, which we give in this definition.

\begin{defn}Let $\mathcal M$ be an $\mathcal L$-structure   and let $A$ be a subset of $\mathcal M$. We define  the \textit{restricted algebraic closure}, denoted by  $racl_{\mathcal M}(A)$, to be   the set of elements $x \in \mathcal M$ such that the orbit
$\{f(x) | f \in Aut(\mathcal M/A)\}$ is finite,  and we define the  \textit{restricted definable closure}, denoted by  $rdcl_{\mathcal M}(A)$, to be the set of elements $x \in \mathcal M$ such that the previous orbit is a singleton; here  $Aut(\mathcal M/A)$ denotes the group of automorphisms of $\mathcal M$ that fix  $A$ pointwise.  To avoid  heaviness of notation,  the subscript $\mathcal M$ will be omitted   if there is no possible confusion. \end{defn}

The following lemma brings together  elementary facts about the  previously defined closures. Its proof  is left to the reader.

\begin{lem}\label{lem-properties-acl} Let $\mathcal M$ be an $\mathcal L$-structure,  and $A,B$  subsets of $\mathcal M$.

$(1)$ $acl(A)$, $dcl(A)$, $acl^\exists(A)$, $dcl^\exists(A)$, $racl(A)$, $rdcl(A)$ are $\mathcal L$-substructures of $\mathcal M$.

$(2)$ $dcl(A) \leq acl(A) \leq racl(A)$, $dcl(A) \leq rdcl(A)$.

$(3)$  $acl(A)=acl^\exists (acl(A))=acl(acl(A))=acl(dcl(A))=dcl(acl(A))=dcl^\exists(acl(A))$.

$(4)$  $A \subseteq B \implies acl(A)\subseteq acl(B)$ and similarly for the other notions of closedness.

$(5)$ If $x \in acl(A)$, then there exists a finite subset $A_0$ of $A$ such that $x \in acl(A_0)$.

$(6)$ If $\mathcal M$ is saturated and $|A|<|\mathcal M|$ then $acl(A)=racl(A)$; similarly for  definable closure. \qed

\end{lem}

Recall that the \textit{type} of a tuple $\bar a \in \mathcal M^n$ over a subset $A$,  denoted $tp(\bar a|A)$,  is the set of formulas $\varphi(\bar x)$ with parameters from $A$ such that  $\mathcal M \models \varphi(\bar a)$, and the \textit{existential type}, denoted $tp^\exists(\bar a|A)$,   is the set of existential formulas $\varphi(\bar x)$ with parameters from $A$ such that  $\mathcal M \models \varphi(\bar a)$. The following proposition is standard, but for completeness we provide a proof of the second property (2) for which we did not find an explicit reference.

\begin{prop} \label{mono-exten} Let $\mathcal M$ be an $\mathcal L$-structure,   $\bar a, \bar b \in \mathcal M^n$ and $A$ a subset of $\mathcal M$.

$(1)$ $tp(\bar a|A)=tp(\bar b|A)$ if and only if there exist an elementary extension $\mathcal N$ of $\mathcal M$ and an automorphism $f \in Aut(\mathcal N/A)$ sending $\bar a$ to $\bar b$.

$(2)$ $tp^\exists(\bar a|A)\subseteq tp^\exists(\bar b|A)$ if and only if there exist an elementary extension $\mathcal N$ of $\mathcal M$ and a monomorphism $f : \mathcal N \rightarrow \mathcal N$, fixing $A$ pointwise and  sending $\bar a$ to $\bar b$.
\end{prop}

\proof  $\;$

$(1)$  See for instance \cite[Theorem 4.1.5]{Marker}.

$(2)$  Clearly, if there is some elementary extension $\mathcal N$ of $\mathcal M$ and  and a monomorphism $f : \mathcal N \rightarrow \mathcal N$, fixing $A$ pointwise and  sending $\bar a$ to $\bar b$, then $tp^\exists(\bar a|A)\subseteq tp^\exists(\bar b|A)$. It remains to show the converse.  Set $\mathcal N_0=\mathcal M$ and let $\mathcal N_1$ be a $|\mathcal M|$-saturated elementary extension of $\mathcal M$.  Using the saturation of $\mathcal N_1$, we get a monomorphism $f_0 : \mathcal N_0 \rightarrow \mathcal N_1$ satisfying $f_0(\bar a)=\bar b$ and fixing $A$ pointwise. Using a similar argument, we build  an elementary chain  $( \mathcal N_i)_{i \in \mathbb N}$,  $\mathcal N_i \preceq \mathcal N_{i+1}$, with  a sequence of monomorphisms  $(f_i : \mathcal N_{i} \rightarrow \mathcal N_{i+1})_{i \in \mathbb N}$ such that  $f_i \restr{\mathcal N_i}=f_{i+1} \restr{\mathcal N_i}$ for every $i \in \mathbb N$. By setting $\mathcal N=\bigcup_{i \in \mathbb N} \mathcal N_i$ and $f=\bigcup_{i \in \mathbb N}f_i$,  we get the required elementary extension and the required monomorphism. \qed

For the reader's convenience, we recall the definition of ultrapowers in the particular case of group theory. An \emph{ultrafilter}  on a set $I$ is a finitely additive
probability measure $\mu : \mathcal {P}(I) \rightarrow \{0,1\}$.
An ultrafilter $\mu$ is called  \emph{nonprincipal} if $\mu(X)=0$
for every finite subset $X \subseteq I$.

Given an ultrafilter $\mu$ on $I$ and a sequence  of groups
$(G_i)_{i \in I }$ we define  an equivalence relation $\sim_{\mu}$
on $\prod_{i \in I}G_i$ by
$$\hat a=(a_i \in G_i)_{ i \in I}\sim_{\mu} \hat b=(b_i \in
G_i)_{ i \in I}~\hbox{ if and only if } \mu(\{i \in
I|a_i=b_i\})=1.
$$

The set of equivalence classes $(\prod_{i \in I}G_i)/\sim_{\mu}$
is endowed with a structure of group by defining
$$
\hat a . \hat b=\hat c ~\hbox{ if and only if } \mu(\{i \in
I|a_i.b_i=c_i\})=1.
$$
The group $(\prod_{i \in I}G_i)/\sim_{\mu}$ is called the
\emph{ultraproduct} of the family $(G_i)_{i \in I }$. When $G_i=G$
for all $i \in I$, $(\prod_{i \in I}G_i)/\sim_{\mu}$ is called an
\emph{ultrapower} and it is denoted simply by $G^*$. If $\mu$ is
nonprincipal, then $G^*$ is called a \emph{nonprincipal
ultrapower}.

\bigskip
\noindent \textbf{Convention.}\textit{  Through this paper we will  consider only ultrapowers on the set of natural numbers; i.e. $I=\Bbb N$ in the previous definition. }
\bigskip

Define $\pi : G \rightarrow {G^*}$ by $\pi (g)=(g_i=g|i \in I)$.
Then $\pi$ is an embedding. Moreover, a theorem of \L os
\cite[Theorem 4.1.9]{Chang-Keisler} claims  that $G$ is an
\emph{elementary subgroup} of $G^*$; that is, any sentence with
parameters from $G$ which is true in $G$ is also true in $G^*$.  In particular, we note that, for any subset $A$ of $G$,  $acl_{G}(A)=acl_{G^*}(A)$ and similarly for  definable closure and their existential correspondents.

Recall that a \textit{countable} model $\mathcal M$ is called \textit{homogeneous} (resp. \textit{$\exists$-homogeneous}), if for any $n \geq 1$, for   any tuples $\bar a$, $\bar b$ of $\mathcal M^n$,  if $tp^{\mathcal M}(\bar a)=tp^{\mathcal M}(\bar b)$ (resp. $tp_{\exists}^{\mathcal M}(\bar a)=tp_{\exists}^{\mathcal M}(\bar b)$) then there exists an automorphism  of $\mathcal M$ which sends $\bar a$ to $\bar b$.  We note,  in particular, that $\exists$-homogeneity implies homogeneity.  For further notions of homogeneity,  we refer the reader to \cite{Hodges(book)93, Marker}.

It is shown in \cite{ould-homogeneity} and \cite{perin-homo} that  nonabelian free groups of finite rank are homogeneous. In the sequel we need the following theorem proved in \cite{ould-homogeneity}. Recall also that a group $G$ is said to be \textit{freely indecomposable} relative to a subgroup $A$, if there is no nontrivial free decomposition  of $G$  such that $A$ is contained in one of the factors.

\begin{thm}\emph{\cite[Proposition  5.9]{ould-homogeneity}} \label{formula-homogeneity} Let $F$ be a nonabelian free group of finite rank  and let $\bar a$  be a tuple of $F$ such that $F$ is freely indecomposable relative to the subgroup generated by $\bar a$. Let $\bar s$ be a basis of $F$. Then there exists a universal formula $\varphi(\bar x)$ such that $F \models \varphi(\bar s)$ and such that for any endomorphism  $f$ of $F$,  if $F \models \varphi(f(\bar s))$ and $f$ fixes $\bar a$ then $f$ is an automorphism.  In particular $(F, \bar a)$ is a prime model of the theory $Th(F, \bar a)$. \qed
\end{thm}

\subsection{Limit groups, modular groups \& abelian JSJ-decompositions}$\;$

Limit groups of  free groups have been introduced by Sela
\cite{Sela-Diophan1} to study  equations over free groups. They can
be seen, geometrically and algebraically,  as \emph{limits}  of
free groups. This class  coincides with the class of \textit{fully
residually-free groups}, a class of groups introduced by Baumslag
\cite{Baum} and studied by Kharlampovich and  Myasnikov
\cite{Kharla-Mias-Irred1, Kharla-Mias-Irred}
and by many other authors.  We start  by giving a  definition which uses ultrafilters in a general context.

\begin{defn} Let $\G$ be a group and $H$ a finitely generated group.  Let $\omega$ be a nonprincipal ultrafilter over $\Bbb N$ and $f =(f_n : H \rightarrow \G)_{n \in \Bbb N}$ a sequence of homomorphisms. Let $\ker_{\omega}(f)$ be the set of elements $h \in H$ such that $\omega(\{n \in \Bbb N| f_n(h)=1\})=1$.  A \textit{$\Gamma$-limit group}   is a group $G$ such that  there exists a finitely generated group $H$, a nonprincipal ultrafilter $\omega$  and a  sequence of homomorphisms $f =(f_n : H \rightarrow \G)_{n \in \Bbb N}$  such that  $G= H/ker_{\omega}(f)$. \qed
\end{defn}

Here is a more standard definition.

\begin{defn} Let $\G$ be a group and $H$ a finitely generated group.  A sequence of homomorphisms $f =(f_n : H \rightarrow \G)_{n \in \Bbb N}$  is called \textit{stable} if, for any $h \in H$, either $f_n(h)=1$ for all but finitely many  $n$, or   $f_n(h) \neq 1$ for all but finitely many $n$.  The \textit{stable
kernel} of $f$, denoted $Ker_{\infty}(f)$,  is the set of elements $h \in H$ such that $f_n(h)=1$ for all but finitely many $n$.  A \textit{$\Gamma$-limit group}  is a group $G$ such that  there exists a finitely generated group $H$   and a stable sequence of homomorphisms $f =(f_n : H \rightarrow \G)_{n \in \Bbb N}$  such that  $G=H/ Ker_{\infty}(f)$. \qed
\end{defn}

The following lemma explains the relation  between the previous notion, which comes  essentially from  geometrical considerations, and the universal theory of the considered group. Its proof can be found  in \cite[Lemma 2.2]{ould-homogeneity} and \cite[Theorem 2.1]{Ould-equa}. For the definition of universal theories, we refer the reader to \cite{Hodges(book)93, Marker} or \cite{Ould-equa} for a quick overview.

\begin{lem}\label{lem-relation-limit-univ}   Let $\G$ be a group and $G$ a finitely generated group. The following properties are equivalent.

$(1)$ $G$ is a  $\Gamma$-limit  group.

$(2)$ $G$ is a  model of the universal theory of $\Gamma$.

$(3)$ $G$ embeds in every nonprincipal ultrapower of $\G$.  \qed
\end{lem}

In dealing with the existential closure in free groups in the next section, we must work with homomorphisms that do not  necessarily fix the subgroup under consideration (in our case $acl^\exists(A)$).  We introduce the following definition which is more appropriate in our context.

\begin{defn} Let $G_1, G_2$  be  groups and $H$ a subgroup of $G_1$. A sequence of homomorphisms $(f_n : G_1 \rightarrow G_2)_{n \in \mathbb N}$  \textit{bounds $H$ in the limit} if for any $h \in H$ there exists a finite subset $B(h)$ of $G_2$ such that $f_n(h) \in B(h)$ for all but finitely many $n$. \qed
\end{defn}
Next theorem  is a slight generalization of similar theorems which appear in  several places \cite{sela-rips-rigid, Sela-Diophan1, Sela-hyp, groves-2007, chloe-these}. As the proof is almost identical, we just give the necessary changes implied by  the previous definition.

Let $\mathcal C$ be a class of subgroups of $G$. By a \textit{$(\mathcal C, H)$-splitting} of $G$ (or a splitting of $G$ over $\mathcal C$ relative to $H$), we understand a tuple $\Lambda=(\mathcal G(V,E), T, \varphi)$, where $\mathcal G(V,E)$ is a graph of groups such that each edge group is in $\mathcal C$ and  $H$ is elliptic, $T$ is a maximal subtree of $\mathcal G(V,E)$ and $\varphi : G \rightarrow \pi(\mathcal G(V,E), T)$ is an isomorphism; here $\pi(\mathcal G(V,E), T)$ denotes the fundamental group of $\mathcal G(V, E)$ relative to $T$. If $\mathcal C$ is the class of abelian groups or cyclic groups, we will just say \textit{abelian splitting} or \textit{cyclic splitting}, respectively.  Splittings of the form $G_1*_CG_2$ or $G_1*_C=\<G, t| c^t=\varphi(c), c\in C\>$ are called \textit{one-edge splittings}.  Given a group  $G$  and a subgroup $H$  of $G$,  $G$ is said to be \textit{freely $H$-decomposable} if $G$ has a nontrivial free decomposition $G=G_1*G_2$ such that $H \leq G_1$.  Otherwise, $G$ is said to be freely $H$-indecomposable.
\begin{thm}\label{relative-splitting} Let $\Gamma$ be a torsion-free hyperbolic group. Let $G$ be a finitely generated group and $H$ a nonabelian subgroup of $G$ such that $G$ is freely $H$-indecomposable. Let $(f_n : G \rightarrow \Gamma)_{n \in \mathbb N}$ be a stable sequence of pairwise distinct homomorphisms with trivial stable kernel and which bounds $H$ in the limit.  Then $G$ admits a nontrivial  abelian splitting  relative to $H$.
\end{thm}

\noindent \textit{Outline of the proof}. Let $S$ be a finite generating set of $\G$ and $(C(\G, S),d)$ the corresponding Cayley graph. Let $D$ be a finite generating set of $G$ and for each $n \in \Bbb N$, define the length $\lambda_n$ of $f_n$ as $\max_{d \in D} |f_n(d)|_S$, where $|.|_S$ denotes the word length relative to $S$. Let $\omega$ be a nonprincipal ultrafilter over $\Bbb N$.  Since the given homomorphisms are pairwise distinct, $\lim_{n \rightarrow \infty} \lambda_n=\infty$. Then $G$ acts on the asymptotic cone $(Con_{\omega}(\G, e, \lambda), d_\omega)$, relative to the sequence of observation points $e=(e_n=1)_{n \in \Bbb N}$, the sequence of scaling factors $\lambda=(\lambda_n)_{n \in \Bbb N}$ and the ultrafilter $\omega$. An argument  similar to the one used in \cite{sela-rips-rigid, chloe-these} shows that the action is superstable, with abelian arc stabilizers and trivial tripod stabilizers. What remains to show in our context is that the action is nontrivial and that $H$ is elliptic.

We claim that $H$ fixes $e$ in $Con_{\omega}(\G, e, \lambda)$.  Since, for any $h \in H$, $\{|f_n(h)|_S|n\in \Bbb N\}$ is bounded, we have
$
d_\omega(e, he)=\lim_\omega \frac{ |f_n(h)|_S}{\lambda_n}=0,
$
and thus $H$ fixes $e$ as claimed.  We claim now that the action is nontrivial. Since $\max_{d \in D} d_\omega(e, de)=1$, $e$ is not a global fixed point. Since $G$ is finitely generated, if  the action is trivial then there is   some global fixed point $e'$, with $e\neq e'$. Then  $H$ will fix the non-degenerate segment $[e,e']$, though it is not abelian; a contradiction with abelianity of arc stabilizers. To get the desired abelian splitting, one may apply \cite{sela-acylin} or  \cite{Guirardel-action}. \qed

\bigskip
The shortening argument is a key tool in Sela's study of limit groups. Roughly speaking, given a sequence of actions of a finitely generated  group $G$ on  the Cayley graph of the torsion-free hyperbolic group $\G$, we get an action of $G$ on some asymptotic cone $C$ of $\G$; by analyzing this action, we can find a particular type of automorphisms, called \textit{modular automorphisms},  of $G$ which shorten the length of the sequence of the actions. Here we briefly recall   modular automorphisms and the shortening argument (in the relative case). For the treatment in the general framework of  hyperbolic groups, we  refer the reader to \cite{weidmann-equa}.

\begin{defn}\label{defi:dehn-twist}
Let $G$ be a  group, and let $\Lambda$ be an abelian  one-edge splitting of $G$ relative to $H$, with edge group $C$. Let $c \in C$. A \textit{Dehn twist} about $c \in C$ is an automorphism $\phi \in Aut(G)$, defined as follows:
\begin{enumerate}
\item if $G = A *_C B$, $H \leq A$, then $\phi(a) = a, \phi(b) = b^c$ for every $a \in A, b \in B$.
\item if $G = A *_C$, $H \leq A$, with stable letter $t$, then $\phi \restr A = id_A$ and $\phi(t) = tc$. \qed
\end{enumerate}
\end{defn}

Let $\Lambda=(\mathcal G(V,E), T, \varphi)$ be a splitting of a group $G$ and $\phi_v$ an automorphism of the vertex group $G_v$, $v \in V$. Suppose that for each $e \in E$ adjacent to $v$, there exists an element $g_e \in G_v$ such that $\phi_v$ restricts to a conjugation by $g_e$ on $G_e$.  Then there exists an automorphism $\phi$ of $G$, called the \textit{standard extension} of $\phi_v$,  which  extends $\phi_v$ (see \cite[Proposition 5.4]{sela-rips-rigid} for more details).

\begin{defn} \label{defi:abelian-modular} Let $\Lambda=(\mathcal G(V,E), T, \varphi)$ be an abelian  splitting of a group $G$ relative to $H$ and $G_v$ an abelian  vertex group. Let $P$ be the subgroup of $G_v$ generated by the incident edge groups.  Any automorphism $\phi_v$ of $G_v$ which fixes $P$ pointwise, and which fixes also  $H$ pointwise,  has a standard extension to $G$. Such an automorphism is called a \textit{modular automorphism of abelian type}. \qed
\end{defn}

Let $\Lambda=(\mathcal G(V,E), T, \varphi)$ be an abelian  splitting of a group $G$ relative to $H$  and  $v \in V$. The vertex $v$ is called of \textit{surface type}, if $G_v$ is isomorphic to the fundamental group of a compact connected surface $S$ with boundary, which is not a disk or a M\"obius band or a cylinder and such that  each edge group $G_e$ incident on $v$ is conjugate to the fundamental group of a boundary component of $S$.

\begin{defn} \label{defi:surface-modular}  Let $\Lambda=(\mathcal G(V,E), T, \varphi)$ be an abelian  splitting of a group $G$ relative to $H$  and  $v \in V$ be a surface type vertex. Any automorphism $\phi_v$ of $G_v$ which restricts to a conjugation by $g_e$ to each incident edge group $G_e$,  and which  fixes  also  $H$ pointwise,  has a standard extension to $G$. Such an automorphism is called a \textit{modular automorphism of surface  type}. \qed
\end{defn}

\begin{defn} Let $G$ be a group and $H$ a subgroup of $G$. The \textit{abelian modular group} of $G$ relative to $H$, denoted $Mod(G/H)$, is the subgroup of $Aut(G/H)$ generated by Dehn twists, modular automorphisms of abelian type and modular automorphisms of surface type. \qed
\end{defn}

We still need a last definition to express the shortening argument:

\begin{defn} Let $G$ be a finitely generated group and $H$ a subgroup of $G$. Let $\G$ be a torsion-free hyperbolic group. Let $B,A$ be finite generating sets of $G,\G$ respectively. A homomorphism $f : G \rightarrow \G$ is said to be \textit{short} relative to $H$ if for any $\sigma \in Mod(G/H)$, one has
$$
\max_{b \in B} |f(b)|_A \leq \max_{b \in B} |f\circ \sigma (b)|_A,
$$
where $|.|_A$ denotes  word length function of $\G$ with respect to $A$. \qed
\end{defn}

\begin{thm}\label{shortening-thm}
Let $\Gamma$ be a torsion-free hyperbolic group with a finite generating set $A$. Let $G$ be a finitely generated group, with a finite generating set $B$,  and $H$ a nonabelian subgroup of $G$ such that $G$ is freely $H$-indecomposable. Let $(f_n : G \rightarrow \Gamma)_{n \in \mathbb N}$ be a stable sequence of pairwise distinct homomorphisms with trivial stable kernel and which bounds $H$ in the limit. Then for any nonprincipal ultrafilter $\omega$, $\omega(\{n \in \Bbb N| f_n\hbox{ is not short}\})=1$.
\end{thm}

\noindent \textit{Outline of the proof}. Let  $(C(\G, A),d)$ be  the corresponding Cayley graph which is hyperbolic. For each $n \in \Bbb N$, let  $\lambda_n=\max_{d \in D} |f_n(d)|_A$. Let $\omega$ be a nonprincipal ultrafilter over $\Bbb N$.  Since the given homomorphisms are pairwise distinct, $\lim_\omega \lambda_n=\infty$. Then $G$ acts on the asymptotic cone $(T, d_\omega)=(Con_{\omega}(\G, e, \lambda), d_\omega)$,  which is a real tree, relative to the sequence of observation points $e=(e_n=1)_{n \in \Bbb N}$, the sequence of scaling factors $\lambda=(\lambda_n)_{n \in \Bbb N}$ and the ultrafilter $\omega$. As in the outline of the proof of Theorem \ref{relative-splitting},  the action is nontrivial, superstable, with abelian arc stabilizers and trivial tripod stabilizers, and $H$ fixes $e$.

By  Rips decomposition (see \cite{bestvina-feighn, sela-acylin} or Guirardel's version \cite{Guirardel-action}), $T$ has a decomposition as a graph of actions $\mathcal A=(\mathcal G(V,E), (T_v)_{v \in V}, (p_e)_{e \in E})$, where each vertex action of $G_v$ is either of symplicial type, or of surface type (IET type) or of abelian type (axial type).

Set $B=\{b_1, \dots, b_q\}$. Let $I$ be the set of indices $i$ such that the segment $[e, b_ie]$ intersects a surface type component, let $J$ be the set of indices $i$ such that $i \not \in I$ and $[e, b_ie]$ intersects an abelian type component; finally, let  $K$ be the set of indices $i$ such that $[e, b_ie]$ lies in a simplicial component.

By using \cite[Proposition 5.2]{sela-rips-rigid},  it is possible to construct a composition of  surface type modular automorphisms $\sigma_1$ such that $d_\omega(e, \sigma_1(b_i)e)<d_\omega(e, b_ie)$ for all $i \in I$ and $\sigma_1(b_i)=b_i$ for all $i \not \in I$.  Let $J' \subseteq I \cup J$ be the set of indices $i$  such that   $[e, \sigma_1( b_i)e]$ intersects an abelian component. In that case, it is possible to find a composition of abelian type modular automorphisms $\sigma_2$ such that $d_\omega(e, \sigma_2 \circ \sigma_1(b_i)e)<d_\omega(e, \sigma_1(b_i)e)$ for all $i \in J'$ and $\sigma_2 \circ \sigma_1(b_i)=\sigma_1(b_i)$ for all $i \not \in J'$.  Finally let $K'$ be the set of indices $i$   such that $[e, \sigma_2\circ \sigma_1(b_i)e]$ intersects a simplicial component. In that case, we cannot ensure the existence of a unique automorphism; however, we show that there exists a subset $U \subseteq \Bbb N$ such that $\omega(U)=1$ and such that for any $n \in U$, there exists a Dehn twist $\tau_n$ such that $d_n(e_n, \tau_n\circ \sigma_2 \circ \sigma_1(f_n(b_i))e_n)<d_n(e_n, \sigma_2\circ \sigma_1(f_n(b_i))e_n)$ for all $i \in K'$ and $d_n(e_n, \tau_n\circ \sigma_2 \circ \sigma_1(f_n(b_i))e_n)=d_n(e_n, \sigma_2\circ \sigma_1(f_n(b_i))e_n)$ for all $i \not \in K'$.

There exists $U_1 \subseteq \Bbb N$ such that  for any $n \in  U_1$, $d_n(e_n, \sigma_1(f_n(b_i))e_n)<d_n(e_n, f_n(b_i)e_n)$ for any $i \in I$ and  $\sigma_1(f_n(b_i))=f_n(b_i)$ for all $i \not \in I$. Similarly, there exists $U_2 \subseteq \Bbb N$ such that  for any $n \in  U_2$, $d_n(e_n, \sigma_2 \circ \sigma_1(f_n(b_i))e_n)<d_n(e_n, \sigma_1(f_n(b_i))e_n)$ for any $i \in J'$ and  $\sigma_2 \circ \sigma_1(f_n(b_i))=\sigma_1(f_n(b_i))$ for all $i \not \in J'$. By taking $\alpha_n=\tau_n \circ \sigma_2 \circ \sigma_1$ and choosing $U'=U \cap U_1 \cap U_2 \subseteq \Bbb N$  we have  $\omega(U')=1$,  and  for any $b_i \in B$, $d_n(e_n, \alpha_n(f_n(b_i))e_n)<d_n(e_n, f_n(b_i)e_n)$ for any $n \in U'$ which proves the desired result.  For more details, the reader can see \cite{wilton-thesis, chloe-these, weidmann-equa, vallino-thesis}. \qed

\bigskip
One of  applications of the shortening argument was the proof by Rips and Sela  \cite{sela-rips-rigid}  of the fact that the modular group has a finite index in the group of automorphisms.  This  can be generalized slightly as follows (see also \cite{chloe-these}).

\begin{thm} \label{finite-index-auto}Let $\G$ be a torsion-free hyperbolic group,  $G$ a finitely generated group,  $H$ a nonabelian subgroup of $G$ such that $G$ is freely $H$-indecomposable. Let   $e : H \rightarrow \G$ be an embedding.  We suppose that there exists at least an embedding of $G$ in $\G$ whose restriction to $H$ is $e$. Then there exists a finite set $\{f_1, \dots,  f_p\}$ of embeddings of $G$ in $\G$,  whose  restriction to $H$ coincides with $e$ and  such that for any embedding $f : G \rightarrow \G$, whose restriction to $H$ coincides with $e$, there exists a modular automorphism $\sigma \in Mod(G/H)$ such that $f \in\{f_1 \circ \sigma, \dots, f_p \circ \sigma\}$.
\end{thm}

\proof    Let $(f_n : G \rightarrow \G)_{n \in \Bbb N}$ be the sequence of all embeddings of $G$ in $\G$ whose restriction to $H$ is $e$. For each $ n \in \Bbb N$, choose a modular automorphism $\sigma_n \in Mod(G/H)$ such that $f_n \circ \sigma_n$ is short. Suppose for a contradiction  that the set $I=\{f_n \circ \sigma_n| n \in \Bbb N\}$ is infinite. Then it is possible to extract a subsequence of pairwise distinct elements from $I$. Clearly such a subsequence is   stable,  has trivial stable kernel and bounds $H$ in the limit. Hence, by Theorem \ref{shortening-thm}  for an infinite set $U \subseteq \Bbb N$, for every $n \in U$, $f_n \circ \sigma_n$ is not short; which is a contraditcion. \qed

\begin{cor} \label{cor-mono} Let $\G$ be a torsion-free hyperbolic group and $H$ a nonabelian subgroup such that $\G$ is freely $H$-indecomposable. Then any monomorphism $f : \G \rightarrow \G$ fixing $H$ pointwise is an automorphism.
\end{cor}

\proof By Theorem \ref{finite-index-auto}, there exists $n,m \in \Bbb N$ such that $n >  m$ and  $f^n=f^m \circ \tau$ for some $\tau \in Mod(\G/H)$.  Therefore $f^{n-m}=\tau$ and thus $f$ is surjective. \qed

\bigskip
One of the important concepts in Sela's study of limit groups is the \textit{shortening quotient}.

\begin{defn} Let $\G$ be a torsion-free hyperbolic group.   Let $G$ be a finitely generated group, $H$ a nonabelian subgroup of $G$ such that $G$ is freely $H$-indecomposable. Let  $f=(f_n : G \rightarrow \G)_{n \in \Bbb N}$ be a stable sequence of pairwise distinct homomorphisms  which bounds $H$ in the limit and such that each $f_n$ is short. The group $SG=G/Ker_{\infty}(f)$ is called a \textit{shortening quotient} of $G$. \qed
\end{defn}

\begin{thm}\label{shortening-quotient} Every shortening quotient is a proper quotient.
\end{thm}

\proof If it is not the case then  the stable kernel is trivial; thus by Theorem \ref{shortening-thm}, for infinitely many $n$, $f_n$ is not short; a contradiction. \qed

Another important application in this context of a more general version of the shortening argument is the proof by Sela \cite{Sela-hyp} of the \textit{descending chain condition} of $\G$-limit groups.

\begin{thm}\label{thm:descending-chain-condition-for-limit-groups} \emph{\cite{Sela-hyp}}
Let $\G$ be a torsion-free hyperbolic group and  $(G_i)_{i \in \Bbb N}$  a sequence of $\Gamma$-limit groups.  If  $(f_i: G_i \rightarrow G_{i+1})_{i \in \Bbb N}$ is a sequence of epimorphisms, then  all but finitely many of them are isomorphisms. \qed
\end{thm}

As it was indicated in the introduction,  a JSJ-decomposition of a group $G$ over a  class of subgroups $\C$ relative to a subgroup $H$ is  a  splitting of $G$ over  $\C$ relative to $H$,  which describes in certain sense all other possible  splittings of $G$ over  $\C$ relative to $H$. Guirardel and Levitt have developed in \cite{Guirardel-Levitt-JSJ1, Guirardel-Levitt-JSJ2} a general framework of JSJ-decompositions that we will use to give the definition and the principal properties.

Given a group $G$ and two $(\C, H)$-splittings $\Lambda_1$ and $\Lambda_2$ of  $G$,  we say that $\Lambda_1$ \textit{dominates} $\Lambda_2$ if every subgroup of $G$ which is elliptic in $\Lambda_1$ is also elliptic in $\Lambda_2$.  A $(\C, H)$-splitting of $G$ is said to be \textit{universally elliptic} if all edge stabilizers in $\Lambda$ are elliptic in any other $(\C, H)$-splitting of $G$.

A \textit{JSJ-decomposition of $G$ over $\C$ relative to  $H$} is an universally elliptic  $(\C, H)$-splitting   dominating all other universally elliptic $(\C,H)$-splittings.  If $\C$ is the class of abelian subgroups, then we  simply say \textit{abelian JSJ-decomposition}; similarly when $\C$ is the class of cyclic subgroups.

It is shown in \cite{Guirardel-Levitt-JSJ1, Guirardel-Levitt-JSJ2} that JSJ-decompositions exist for finitely presented groups. Here we will use  existence and properties  of JSJ-decompositions in the framework of finitely generated torsion-free CSA-groups proved in  \cite{Guirardel-Levitt-JSJ2}.

Given a surface $\Sigma$, a \textit{boundary subgroup} of the fundamental group $\pi_1(\Sigma)$ is a subgroup conjugate to the fundamental group of a boundary component. An \textit{extended boundary subgroup} of $\pi_1(\Sigma)$ is a subgroup of a boundary subgroup.

Let $G$ be a group and  $\Lambda$ a $(\C, H)$-splitting of $G$. A vertex stabilizer $G_v$ in $\Lambda$ is called of \textit{QH surface type} if it is isomorphic to the fundamental group $\pi_1(\Sigma)$ of a surface $\Sigma$ such that images of incident edge groups are extended boundary subgroups and every conjugate of $H$ intersects $G_v$ in an extended boundary subgroup.  A boundary component $C$ of $\Sigma$ is \textit{used} if there exists an incident edge group, or a subgroup of $G_v$ conjugate to $H$ whose image in $\pi_1(\Sigma)$ is contained with finite index in $\pi_1(\Sigma)$.

A vertex stabilizer $G_v$ in $\Lambda$ is said to be \textit{rigid} if it is elliptic in every $(\C, H)$-splitting of $G$. Otherwise it  is called \textit{flexible}.

Recall that a group  is called \textit{CSA} if every maximal abelian subgroup is malnormal.  It is a general fact  that if $\G$ is a torsion-free hyperbolic group then $\G$-limit groups  are torsion-free and CSA. The following theorem is an application of results  of  \cite{Guirardel-Levitt-JSJ2} in our particular context.

\begin{thm} \label{JSJ-decomp} \emph{\cite[Theorem 11.1]{Guirardel-Levitt-JSJ2}}  Let $G$ be a torsion-free finitely generated   CSA-group and $H$ a subgroup of $G$ such that $G$ is $H$-freely indecomposable. Then abelian JSJ-decompositions of $G$ relative to $H$ exist and their nonabelian flexible vertices are of QH surface type with every boundary component used. \qed
\end{thm}

Since boundary subgroups are cyclic, it follows that if $H$ is nonabelian then $H$ is contained in a conjugate of a rigid group in any abelian JSJ-decomposition of $G$ relative to $H$.  Hence, without loss of generality,  in the rest of this paper we may assume  that JSJ-decompositions used by us have the property that $H$ is contained in a rigid vertex group.  Since, we will use only properties that are satisfied by all  JSJ-decompositions, by  misuse of language we will use the term  \textit{the JSJ-decomposition} rather than a JSJ-decomposition. Through this paper we will use  the following two  simple properties of JSJ-decompositions.

\begin{lem} \label{mod-preservation} Let $G$ be a finitely generated torsion-free CSA-group and $H$ a nonabelian subgroup of $G$ such that $G$ is $H$-freely indecomposable. Let $\Lambda$ be the abelian JSJ-decomposition of $G$ relative to $H$. Then any  automorphism  from $Mod(G/H)$ fixes pointwise the vertex group containing $H$ in $\Lambda$.
\end{lem}

\proof Let $G(H)$ be the vertex group of $\Lambda$ containing $H$. Since $G(H)$ is rigid it is elliptic in any abelian splitting of $G$ relative to $H$. Let $\sigma \in Mod(G/H)$. Suppose that $\sigma$ is a Dehn twist and let $G=G_1*_CG_2$ or $G=L*_C$ be the corresponding one-edge abelian splitting.  Since $H \leq G_1$ or $H \leq L$ and $H \leq G(H)$ which is elliptic, it follows that $G(H) \leq G_1$ or $G(H) \leq L$ which is the desired conclusion.  Using a similar argument, if $\sigma$ is an automorphism of surface type or abelian type then it fixes $G(H)$ pointwise. \qed

\begin{lem} \label{mod-quot}Let $G$ be a finitely generated torsion-free CSA-group and $H$ a nonabelian subgroup of $G$ such that $G$ is $H$-freely indecomposable. Let $f=(f_n : G \rightarrow \Gamma)_{n \in \mathbb N}$ be a stable sequence of pairwise distinct homomorphisms  with trivial stable kernel and which bounds $H$ in the limit. For each $n \in \Bbb N$ choose $\sigma_n \in Mod(G/H)$ such that $f_n \circ \sigma_n$ is short. Let $SG$ be the corresponding  shortening quotient  and $\pi : G \rightarrow SG$  the natural map.  Then the restriction of $\pi$ to the vertex group $G(H)$ containing $H$  in the abelian JSJ-decomposition  of $G$ relative to $H$ is injective.
\end{lem}

\proof By Lemma \ref{mod-preservation}, for every $g \in G(H)$, $f_n \circ \sigma_n(g)=f_n(g)$ and the required conclusion follows.  \qed

\bigskip
All the previous properties of  JSJ-decompositions are widely sufficient in  our context of $\G$-limit groups. However for torsion-free hyperbolic groups themselves,  we  need some additional properties. Let $G$ be a group and  $\Lambda$ a $(\C, H)$-splitting of $G$. We say that a boundary subgroup $B$ of a surface type vertex group $G_v$ is  \textit{fully used}  if there exists an incident edge group, or a subgroup of $G_v$ conjugate to $H$, which coincides with $B$.  

Let  $\Lambda$  be an abelian splitting of $G$ (relative to $H$) and  $G_v$ be a vertex group of $\Lambda$.  The \textit{elliptic abelian neighborhood}  of $G_v$ is the subgroup generated by the elliptic elements that commute with nontrivial elements of $G_v$.  By  \cite[Proposition 4.26]{gui-limit}  if $G$ is commutative transitive then any  abelian splitting $\Lambda$  of $G$ (relative to $H$) can be transformed to an abelian splitting $\Lambda'$ of $G$ such that the underlying graph is the same as that of $\Lambda$ and for any vertex $v$, the corresponding new vertex group $\hat G_v$ in $\Lambda'$ is the elliptic abelian neighborhood of $G_v$ (similarly for edges). In particular any edge group of $\Lambda'$ is malnormal in the adjacent vertex groups  and any boundary subgroup of a surface type vertex group is fully used. We call that transformation the \textit{malnormalization} of $\Lambda$.  If $\Lambda$ is a (cyclic or abelian) JSJ-decomposition of $G$ and $G$ is commutative transitive then the malnormalization of $\Lambda$ will be called a \textit{malnormal} JSJ-decomposition. If $G_v$ is a rigid vertex group then we call $\hat G_v$ also rigid; similarly for abelian and surface type vertex groups. Strictly speaking a malnormal JSJ-decomposition  is not a JSJ-decomposition in the sense  of \cite{Guirardel-Levitt-JSJ1}, however it shares the most important properties with  JSJ-decompositions that we need. Hence we get the following which summarizes several properties    sufficient for our purpose.

\begin{thm} \label{thm-JSJ-hyp}Let $G$ be a torsion-free finitely generated CSA-group and $H$ a nonabelian subgroup of $G$ such that $G$ is $H$-freely indecomposable. Then malnormal  abelian JSJ-decompositions of $G$ relative to $H$ exist and satisfy the following properties.

$(1)$  Flexible vertices are of  QH surface type with every boundary component fully used.

$(2)$ Every edge group is maximal abelian in its endpoints vertex groups.

$(3)$ $H$ is contained in a rigid vertex group. \qed
\end{thm}

\bigskip
We end with the definition of   \textit{generalized JSJ-decomposition}. First, split  $\G$ as a free product $\G=\G_1*\G_2$, where $H \leq \G_1$ and $\G_1$ is freely $H$-indecomposable (relative Grushko-Kurosh decomposition). Then, define the generalized (cyclic) JSJ-decomposition of $\G$ relative to $H$ as the (cyclic) splitting obtained by adding $\G_2$ as a new vertex group to the (cyclic) JSJ-decomposition of $\G_1$ (relative to $H$). The notion of a generalized malnormal (cyclic) JSJ-decomposition is defined in a similar way.

\bigskip
Recall that a group is said to be equationally noetherian if any system of equations in finitely many variables is equivalent to a finite subsytem. For more details on this notion, we refer the reader to \cite{Baum-Rem-algebraic}. A theorem of Sela  \cite[Theorem 1.22]{Sela-hyp} states that any system of equations without
parameters  in finitely many variables is equivalent in a torsion-free
hyperbolic group to a finite subsystem. The previous property is equivalent,
when the group under consideration $G$ is finitely generated, to the fact
that $G$ is equationally noetherian (for more details see the end of section 2 in \cite{ould-homogeneity}). Hence a torsion-free hyperbolic group is equationally noetherian. This was generalized by C. Reinfeldt and R. Weidmann   \cite{weidmann-equa}  to general hyperbolic groups.

\begin{thm} \emph{\cite{weidmann-equa}} \label{weidmann-equa} A  hyperbolic group is equationally noetherian. \qed
\end{thm}

\section{Constructibility from the algebraic closure}

As noticed  before, if $A$ is a subset of $\G$ then $acl(A)$ and $acl(\<A\>)$ coincide, similarly with the other notions of closures,   thus without loss of generality we may assume that $A$ is always a subgroup. First we treat the case of abelian subgroups.  We denote by $C_G(A)$ the centralizer of $A$ in $G$.

\begin{lem} \label{acl-abelian}Let $G$ be a torsion-free CSA group whose abelian subgroups are cyclic. Let $A$ be a nontrivial abelian subgroup of $G$. Then $racl(A)=acl(A)=acl^\exists(A)=dcl^\exists(A)=dcl(A)=rdcl(A)=C_G(A)$.
\end{lem}

\proof   We first show that $racl(A) \leq C_G(A)$. Let $g \in racl(A)$, $a \in A$, $g \neq 1$, $a \neq 1$.  Let $\pi_n$  be the conjugation by $a^n$, $n \in \Bbb N$. Hence the set $\{\pi_n(g)| n \in \Bbb N\}$ is finite. Thus $[a^{n-m}, g]=1$ for some $n, m \in \Bbb N$, $n \neq m$ . Since $G$ is torsion-free and  CSA, commutativity is a transitive relation on the set of nontrivial elements, thus $[g, a]=1$. Therefore $g \in C_G(A)$ as required.

Now we show that $C_G(A) \leq dcl^\exists(A)$.  Since  $C_G(A)$ is cyclic, there exists $b \in G$ such that $C_G(A)=\<b\>$.  Let  $a \in A$, $a \neq 1$ and  $m \in \Bbb Z$ such that $b^m=a$. Therefore $b$ satisfies the equation $x^m=a$.  Since $G$ is torsion-free and commutative transitive, $b$ is the unique element satisfying  $x^m=a$. Hence $b \in dcl^\exists(A)$ and thus $C_G(A) \leq dcl^\exists(A)$ as required.  We conclude by the inclusions given by Lemma \ref{lem-properties-acl}. \qed

\bigskip

Since torsion-free hyperbolic groups are CSA, the previous lemma holds for them. Also note that if $G$ is nonabelian  then the algebraic closure of the trivial element is trivial. Indeed by taking $a,b \in G$ with $[a,b]\neq 1$ we have $acl(1) \leq acl(\<a\>) \cap acl(\<b\>)=1$.

Recall that an $\mathcal L$-subtructure $\mathcal N$ of an $\mathcal L$-structure $\mathcal M$ is said to be \textit{existentially closed}, abreviated e.c.,  if for any  existential  formula $\varphi$ with parameters from $\mathcal N$, if $\mathcal M \models \varphi$,  then $\mathcal N \models  \varphi$.   To avoid repeating some proofs, we introduce the following weak notion of existential closedness, of independent interest. A subset $A$ of an $\mathcal L$-structure is said to be \textit{finitely existentially closed} if $acl^\exists(A)=A$.  For instance a nontrivial centralizer in a torsion-free hyperbolic group is  finitely existentially closed (Lemma \ref{acl-abelian} above). It follows immediately that a finitely existentially closed subset is in fact an $\mathcal L$-substructure, so in in the particular context of groups it  is a subgroup.  The first aim of this section is a proof of next  theorem.  First we give a definition.

\begin{defn} Let $G$ be a group,  $A$ a subgroup and $\mathcal C$ a class of subgroups. By induction on $n$, define $\mathcal D_0=\{A\}$, $\mathcal D_{n+1}=\mathcal D_n \cup \{B_1*_CB_2, B*_C| B_1, B \in \mathcal D_n, B_2 \leq G, C \in \mathcal C\}$. We say that $G$ is \textit{constructible from $A$ over $\mathcal C$}, if there exists $n \in \mathbb N$ such that $G \in \mathcal D_n$. \qed
\end{defn}

\begin{thm}\label{thme} Let $\G$ be a torsion-free hyperbolic group and $A$ a nonabelian finitely existentially closed subgroup of $\G$. Then $\G$ is constructible  from $A$  over cyclic subgroups.  In particular $A$  is    finitely generated, quasiconvex (and hyperbolic).
\end{thm}

Since for any subset $A$,   $acl(A)$ is finitely existentially closed (Lemma \ref{lem-properties-acl}(3)), Theorem \ref{thme} implies Theorem \ref{thm-princip-1}.  It is shown in \cite{chloe-these} that, given a torsion-free hyperbolic group $\G$, if $A$ is an elementary subgroup then $\G$ has a structure of a hyperbolic tower over $A$ and in particular $A$ is finitely generated, quasiconvex and hyperbolic. Theorem \ref{thme} allows to deduce these last properties  which generalize to existentially closed subgroups, too. Indeed,  since an existentially closed subgroup is in particular finitely existentially closed, we obtain the following.

\begin{cor} An existentially closed subgroup of a torsion-free hyperbolic group is finitely generated, quasiconvex (and hyperbolic). \qed
\end{cor}

The first part  of this section is devoted to the proof of Theorem \ref{thme}. We start with the following lemma of general interest.

\begin{lem} \label{lem-equa-hom} Let $G$ be an  equationally noetherian group. Let $G^*$ be an elementary extension of $G$. Let $P$ be a subset of $G$. Let $K$ be a finitely generated subgroup of $G^*$ such that $P \subseteq K$. Then there exists a finite subset $P_0 \subseteq P$ such that for any homomorphism $f: K \rightarrow G^*$,  if $f$ fixes  $P_0$ pointwise then $f$ fixes $P$ pointwise.
\end{lem}

\proof   Let $\bar g$ be a generating tuple of $K$. Write $P=\{p_i | i \in \mathbb N\}$.  Then for every $i \in \mathbb N$,  there exists a word $w_i(\bar x)$ such that $p_i=w_i(\bar g)$.  Since $G$ is equationally noetherian and $P \subseteq G$, there exists $n \in \mathbb N$ such that
$$
G^*\models \forall \bar x (( p_0=w_0(\bar x) \wedge \dots \wedge p_n=w_n(\bar x) ) \implies p_i=w_i(\bar x)), \leqno (1)
$$
for any $i \in \mathbb N$.

Let $P_0=\{p_0, \dots, p_n\}$ and let $f: K \rightarrow G^*$ be a homomorphism  such that $f(p_i)=p_i$ for every $0 \leq i \leq n$. Therefore $p_i=f(p_i)=w_i(f(\bar g))$ for any $0 \leq i \leq n$. Hence, by $(1)$, $p_i=w_i(f(\bar g))$  for any $i \in \mathbb N$, thus $p_i=f(p_i)$ for any $i \in \mathbb N$, as required. \qed

\begin{prop} \label{prop1} Let $\Gamma$ be a torsion-free hyperbolic group and $A$ a nonabelian finitely existentially closed subgroup of $\G$.  Let $\G^*$ be a nonprincipal ultrapower of  $\Gamma$. Let $K \leq {\G^*}$ be a finitely generated subgroup such that $A \leq K$ and such that $K$ is $A$-freely indecomposable. Then one of the following cases holds.

\smallskip
$(1)$ Let $\Lambda$ be the abelian JSJ-decomposition of $K$ relative to $A$. Then  the vertex group containing $A$ in $\Lambda$ is exactly $A$.

$(2)$ There exists a finitely generated subgroup $L \leq {\G^*}$ such that $A \leq L$ and a non-injective epimorphism  $f : K \rightarrow L$ satisfying:

$\; (2)(i)$  $f$ sends $A$ to $A$ pointwise;

$\; (2)(ii)$  if $\Lambda$ is the abelian JSJ-decomposition of $K$ relative to $A$,  then the restriction of $f$ to the vertex group containing $A$ in $\Lambda$ is injective.
\end{prop}

\proof

Let $\bar d=(d_1, \dots, d_p)$ be a finite generating tuple of $K$. Let $\bar x=(x_1, \dots, x_p)$ be a new tuple of variables and set
$$
S(\bar x)=\{w(\bar x)| K \models w(\bar d)=1\},
$$
where $w(\bar x)$ denotes a word on $\bar x$ and their inverses.

Since $\G$ is equationally noetherian and $\G^*$ is an elementary extension of $\G$,  there exist words $w_1(\bar x), \dots, w_m(\bar x)$  from $S(\bar x)$,  such that
$$
\G^* \models \forall  \bar x (w_1(\bar x)=1 \wedge \dots \wedge w_m(\bar x)=1 \implies w(\bar x)=1),
$$
for any $w \in S(\bar x)$.

By Lemma \ref{lem-equa-hom}, there exists a finite subset $P_0=\{p_1, \dots, p_q\} \subseteq A$,   such that for any homomorphism   $f: K \rightarrow \G$,  if $f$ fixes $P_0$ pointwise then $f$ fixes $A$ pointwise.  Let $p_1(\bar x), \dots, p_q(\bar x)$ be words such that $p_i(\bar d)=p_i$ for every $1 \leq i \leq q$. Set
$$
\phi(\bar x) := w_1(\bar x)=1 \wedge \dots \wedge w_m(\bar x)=1\wedge p_1(\bar x)=p_1 \wedge \dots \wedge p_q(\bar x)=p_q.
$$

We conclude that  any map $f : K \rightarrow \G$ satisfying  $\G \models \phi(f(\bar d))$  extends to a homomorphism which  fixes $A$ pointwise, that we still denote $f$.

Let $(v_i(\bar x)| i \in \mathbb N)$ be the list of reduced words such that $K \models v_i(\bar d) \neq 1$. For  $m \in \mathbb N$,  we set
$$
\varphi_m(\bar x):= \phi(\bar x) \wedge \bigwedge_{0 \leq i \leq m}v_i(\bar x) \neq 1.  \leqno (*)
$$

Suppose first  that there exists $m \in \mathbb N$,  such that for  any map  $f: K \rightarrow \Gamma$  for which $\G \models \varphi_m(f(\bar d))$,   $f$ is an embedding. We claim that,  in that case, the vertex group $B$ containing $A$ in the abelian JSJ-decomposition of $K$ relative to $A$ is exactly $A$. Thus we obtain conclusion $(1)$ of the proposition.

Let $\bar b$ be a finite generating tuple of $B$.  Then there exists a tuple of words $\bar w(\bar x)$ such that $\bar b =\bar w(\bar d)$. We claim that the formula
$$
\psi(\bar y):=\exists \bar x (\varphi_m(\bar x)\wedge \bar y =\bar w(\bar x)),
$$
has only finitely many realizations in $\Gamma$.

Let $\bar c $ in $\G$ such that $\G \models \varphi(\bar c)$. Hence there exists an embedding $f: K \rightarrow \G$, fixing pointwise $A$,  such that $\bar c= \bar w(f(\bar d))$.  Thus the subgroup generated by $\bar c$ is the image of $B$ by $f$.

By  Theorem \ref{finite-index-auto}, there exist finitely many embeddings $h_1, \dots, h_k$, fixing $A$ pointwise,  such that for any embedding $h : K \rightarrow \G$, there exists a modular automorphism $\tau \in Mod(\G/A)$ such that $h\circ \tau =h_i$.  Since any modular automorphism fixes $B$ pointwise (Lemma \ref{mod-preservation}),  we find  $\bar c= f(\bar b) \in  \{h_1(\bar b), \dots, h_k(\bar b)\}$,  thus we get the required conclusion.  Since $\G^* \models \varphi(\bar b)$,  we conclude that $B \leq acl^\exists(A)=A$ as claimed.

Suppose now that for every  $m \in \mathbb N$, there exists a non-injective homomorphism $f: K \rightarrow \G$ such that $\G \models \varphi_m(f(\bar d))$. Therefore, we get a stable sequence $(f_m : K \rightarrow \Gamma)_{m \in \Bbb N}$ of pairwise distinct homomorphisms  with trivial stable  kernel.

For each $n \in \mathbb N$,  choose a modular automorphism $\tau_n \in Mod(K|A)$ such that $h_n=f_n \circ \tau_n$ is short relative to $A$. Hence, we extract  a stable subsequence $(h_m : K \rightarrow \Gamma)_{m \in \Bbb N}$ of pairwise distinct homomorphisms. Let $L$ be the corresponding shortening quotient, which is embeddable in $^*\Gamma$ and contains $A$ and let $f : K \rightarrow L$ be the quotient map.  By Theorem \ref{shortening-quotient}  $L$ is a proper quotient.  We see also that $f$ sends $A$ to $A$ pointwise. Since the stable kernel of $(f_n : K \rightarrow \Gamma)$ is trivial and since every modular automorphism fixes $B$ pointwise, the restriction of $f$ to $B$ is injective (Lemma  \ref{mod-quot}). Hence we obtain conclusion $(2)$ of the proposition. This ends the proof of the proposition.  \qed


\begin{cor} \label{cor-const}Let $\G$ be a torsion-free hyperbolic group and $A$ a nonabelian finitely existentially closed subgroup of $\G$.  Let $\Gamma^*$ be a nonprincipal ultrapower of  $\G$. Let $K \leq {\G^*}$ be a finitely generated subgroup containing $A$. Then $K$ is constructible  from $A$ over  abelian subgroups.
\end{cor}

\proof We construct  a sequence  $K=K_0, K_1, \dots, K_n$ of finitely generated subgroups of $\G^*$, with epimorphisms $f_i: K_i \rightarrow K_{i+1}$ satisfying:

\smallskip
$(i)$ $f_i$ sends $A$ to $A$ pointwise,

$(ii)$ either $K_{i+1}$ is a free factor of $K_i$ and $f_i$ is just the retraction that kills the complement,  or   the restriction of $f$   to the vertex group containing $A$ in the abelian JSJ-decomposition of $K_i$ relative to $A$ is injective,

$(iii)$ if $\Lambda$ is the abelian JSJ-decomposition of $K_n$,  then the vertex group containing $A$ in $\Lambda$ is exactly $A$.

\smallskip

 We put $K_0=K$. Suppose that $K_i$ is constructed. If $K_i$ is freely decomposable relative to $A$,  then we set $K_i=K_{i+1}*H$ with $A \leq K_{i+1}$ and $K_{i+1}$ freely $A$-indecomposable. We define $f_i :K_i \rightarrow K_{i+1}$ to be the retraction that kills $H$.

If $K_i$ is freely $A$-indecomposable,  then one of the cases of Proposition \ref{prop1} is fulfilled. If $(1)$ holds,  then this terminates the construction of the sequence. Otherwise, $(2)$ of Proposition \ref{prop1} holds and we get   $K_{i+1}\leq \G^*$ and $f_{i} : K_i \rightarrow K_{i+1}$ satisfying  $(2)(i)$-$(ii)$ of Proposition \ref{prop1}.

Using the descending chain condition  on $\Gamma$-limit groups (Theorem \ref{thm:descending-chain-condition-for-limit-groups}), the sequence terminates. Let $K_n$ be the last element in the sequence.  Hence, property $(iii)$ is satisfied. We show by inverse induction on $i$, that $K_i$ satisfies the conclusion of the corollary. Since $A$ is exactly the vertex group containing $A$ in the abelian JSJ-decomposition of $K_n$ relative to $A$, it follows that $K_n$ can be constructed from $A$ by a sequence of amalgamated free products and HNN-extensions along abelian subgroups. Hence $K_n$ satisfies the conclusion of the corollary.

 Suppose that $K_{i+1}$ satisfies the conclusions of the corollary. By construction, either  $K_i=K_{i+1}*H$, in which case $K_i$ satisfies the conclusion of the corollary, or the restriction $f_i$   to the vertex group $V$ containing $acl(A)$ in the abelian JSJ-decomposition of $K_i$ relative to $A$ is injective.  By induction, $K_i$ satisfies the conclusions of the corollary. Since $f_i(V)$ contains $A$ and $f_i$ sends $A$ to $A$ pointwise,  $f_i(V)$ is constructible from $A$ by a sequence of amalgamated free products and HNN-extensions along abelian subgroups.  Since the restriction of $f_i$ to $V$ is injective, it follows that $V$ itself is constructible  from $A$ by a sequence of free products and HNN-extensions along abelian subgroups. Therefore $K_{i}$ satisfies the conclusion of the corollary.  Hence $K$ is constructible  from $A$ by a sequence of amalgamated free products and HNN-extensions along abelian subgroups; thus the corollary is proved. \qed

 \bigskip
 Following \cite{ventura-martino-compress}, a subgroup $A$ of a free group $F$ is \textit{compressed} if whenever $A \leq K$,  with $K$ finitely generated,  then $rk(A) \leq rk(K)$; here $rk(H)$ denotes the rank of $H$.

 \begin{cor} \label{acl-compress} Let $F$ be a free group of finite rank and $A$ a nonabelian subgroup of $F$. Then $acl(A)$ is compressed.
  \end{cor}

  \proof By Corollary \ref{cor-const} if $acl(A) \leq K$, with $K$ finitely generated,  then $K$ is constructible from $acl(A)$ over cyclic subgroups. Let $K=B_1*_CB_2$ with $acl(A) \leq B_1$ and $C=\<c\>$. By \cite[Theorem 1.1]{ould-not-cyc}, $c$ is either primitive in $B_1$ or $B_2$.  Therefore $rk(B_i) \leq rk(K)$ for $i=1,2$. Similarly, if $K=B*_C$ then $rk (B)\leq rk(K)$; a consequence of \cite[Theorem 1.1]{ould-not-cyc}. Hence, by induction we get that  $rk(acl(A))\leq rk(K)$.  \qed

\bigskip
\noindent \textbf{Proof of Theorem \ref{thme}.} The fact that $\G$ is constructible from $A$ over cyclic subgroups follows from Corollary \ref{cor-const}. Since $\G$ is finitely generated, any vertex group in any cyclic splitting of $\G$ is finitely generated. Thus by induction and using the fact that $\G$ is constructible from $A$ over cyclic subgroups we find that $A$ is finitely generated. The same argument combined with the following theorem shows that $A$ is quasiconvex and in particular hyperbolic.  \qed

\begin{thm}\emph{\cite[Proposition 4.5]{kapovich-two} } Let $\G$ be a hyperbolic group. Suppose that $\Lambda$ is a cyclic splitting of $\G$ with a finite underlying graph. Then all vertex groups
of $\Lambda$ are quasiconvex in $\G$ and word-hyperbolic themselves. \qed
\end{thm}

Note that in general $acl^\exists(A)$ is not finitely existentially closed, thus Theorem \ref{thme} cannot be applied to  existential algebraic closure. The rest of this section is devoted to show that free groups of finite rank are constructible from the existential algebraic closure. The general case of torsion-free hyperbolic groups is studied in \cite{vallino-thesis}.


\begin{thm} \label{thm-const-exis} Let $F$ be a free group of finite rank and $A$ a nonabelian subgroup of $\G$.  Let $K$ be a finitely generated subgroup of $F$ containing $acl^\exists(A)$. Then $K$ is constructible from $acl^\exists(A)$ over cyclic subgroups.
\end{thm}

First we prove the following general key proposition of independent interest.

\begin{prop} \label{prop-existent-const} Let $G$ be a finitely generated equationally noetherian group and let $A$ be a subgroup of $G$. Let $K \leq G$ be finitely generated and suppose that  $acl^\exists(A)$ is a proper subgroup of $K$. Then there exists a stable sequence of pairwise distinct homomorphisms $({h_n: K \rightarrow F})_{n \in \mathbb N}$ with trivial stable kernel and which bounds $acl^\exists(A)$ in the limit.
\end{prop}

In what follows we fix a finitely generated equationally noetherian group $G$ and $A$ a subgroup of $G$. We fix a finite generating set of $G$ and we denote by $B_r$ the ball of radius $r$ with respect to the word distance induced by the fixed generating set. We denote by $Mon(G/A)$ the monoid of monomorphisms of $G$ fixing $A$ pointwise.  We introduce  the following definition.
\begin{defn}\label{defi:strongly-converges}  Let $G^*$ be an elementary extension of $G$ and let $C$ be a finitely generated subgroup of $G^*$. A stable sequence $(f_n: C \rightarrow G)_{n \in \mathbb N}$ with trivial stable kernel \textit{strongly converges} to $C$ if it satisfies the following properties:
\begin{enumerate}
 \item\label{item:fixes-who-is-inside} for any $c \in C \cap G$, $f_n(c)=c$ for all but finitely many $n$;
 \item\label{item:moves-who-is-outside} for any $c \in C$, for any $b \in G$, if $f_{n_k}(g)=b$ for some subsequence  $(n_k)_{k \in \mathbb N}$, then $g=b$.
\end{enumerate}
\end{defn}

\begin{lem} \label{thm:existence-of-converging-sequence} Let $G^*$ be an elementary extension of $G$ and let $C \leq G^*$ be finitely generated. Then there exists a stable sequence of homomorphisms  $(f_n: C \rightarrow G)_{n \in \mathbb N}$ strongly converging to $C$.
\end{lem}

\textit{Proof. }Let
$$C=\sq{c_1,\dots, c_t}{w_i(\bar c)=1, i \in \mathbb N}$$
be a presentation of $C$. Since $G$ is equationally noetherian, there exists a finite number of words $w_0, \dots, w_p$ such that
$$G \models \forall \bar x((w_0(\bar x)=1 \wedge \dots \wedge w_p(\bar x)=1) \Rightarrow w_i(\bar x)=1)$$
for any $i \in \mathbb N$.

Enumerate the following sets:
$$G\smallsetminus \{1\}=(a_i)_{i \in \mathbb N},$$
$$(G \cap C)\smallsetminus \{1\}=(b_i)_{i \in \mathbb N}=(b_i(\bar c))_{i \in \mathbb N}$$
and
$$C \smallsetminus  \{1\}=(v_i(\bar c))_{i \in \mathbb N},$$
and
$$C \smallsetminus G=(d_i(\bar c))_{i \in \mathbb N}.$$

Since $G$ is an elementary subgroup of $G^*$, for any $n \geq 0$ there exists $\bar c_n$ in $G$ such that
\begin{equation}\label{eqn:stable-sequence}
G \models  \bigwedge_{0\leq i \leq p}w_i(\bar c_n )=1 \wedge  \bigwedge_{0\leq i \leq n}v_i(\bar c_n) \neq 1
\end{equation}
and
\begin{equation}\label{eqn:approximation}
G \models  \bigwedge_{0\leq i \leq n}b_i=b_i(\bar c_n) \wedge \bigwedge_{0 \leq i \leq n, 0 \leq j \leq n} d_i(\bar c_n) \neq a_j.
\end{equation}

Define  $f_n(\bar c)=\bar c_n$ and we show that the sequence $(f_n)_{n \in \mathbb N}$ satisfies properties \ref{item:fixes-who-is-inside} and \ref{item:moves-who-is-outside} of Definition \ref{defi:strongly-converges}.

The sequence $(f_n)_{n \in \mathbb N}$ is stable and has a trivial stable kernel by equation \eqref{eqn:stable-sequence}. Let $g \in C \cap G$. Then there exists $m$ such that $g=b_m=b_m(\bar c)$.  By equation \eqref{eqn:approximation}, we have $f_n(b_m(\bar c))=b_m(\bar c_n)=b_m$ for any $n \geq m$; thus $f_n(g)=g$ for all but finitely many $n$, so we have property \ref{item:fixes-who-is-inside}.

Now, let $g \in C$ and $b \in G$ such that there exists a subsequence $(n_k)_{k \in \mathbb N}$ with $f_{n_k}(g)=b$ for any $k \geq 0$.  Let $s$ be such that $b=a_s$. Suppose for a contradiction that $g \not \in G$. Then there exists $r$ such that $g=d_r(\bar c)$.  Let $n \geq \max \{r,s\}$. By equation \eqref{eqn:approximation}, we have $f_n(g)=f_n(d_r(\bar c))=d_r(\bar c_n)\neq a_s$. Therefore for $n_k$ large enough we have $f_{n_k}(g)\neq b$; a contradiction.

Hence $g \in G$ and  in particular $g \in C \cap G$.  By property \ref{item:fixes-who-is-inside} we get $f_n(g)=g$ for all but finitely many  $n$ and in particular $g=b$ as required, so property \ref{item:moves-who-is-outside} is proved. \qed

\begin{lem}  \label{thm:equivalence-acl-ex-extension-converging} The following properties are equivalent for any finite subset $C \subseteq G$:
\begin{enumerate}
\item\label{item:subset-acl-ex} $C \subseteq acl^{\exists}(A)$;
 \item\label{item:subset-extension-ex} there exists a finite subset $B(C) \subseteq G$ such that for any elementary extension $G^*$ of $G$ and for any $f \in Mon(G^*/A)$, $f(C) \subseteq B(C)$;
 \item\label{item:converging-ex} there exists $r > 0$ such that for any elementary extension $G^*$ of $G$, for any $f \in Mon(G^*/A)$,  for any sequence $(g_n: f(G) \rightarrow G)_{n \in \mathbb N}$ which strongly converges to $f(G)$, $(g_n\circ f)(C) \subseteq B_r$ for all but finitely many $n$.
\end{enumerate}
\end{lem}

\textit{Proof. } $(1) \Rightarrow (2)$. This follows immediately from the definition of $acl^{\exists}(A)$.

$(2) \Rightarrow (3)$. Let $B(C)$ be the given subset. Let
$$r= \max \{|g|\;; \; g \in B(C)\},$$
where $|.|$ is the word length with respect to the finite generating set of $G$.
Let $G \preceq G^*$, let $f \in Mon(G^*/A)$ and let $(g_n: f(G) \rightarrow G)_{n \in \mathbb N}$ be a sequence strongly converging to $f(G)$.  Let $c \in C$. Hence $f(c)=b \in B(C) \subseteq G$ and $b \in G \cap f(G)$.  Since $(g_n)_{n \in \mathbb N}$ strongly converges to $f(G)$, we have $g_n(b)=b$ for all but finitely many $n$. Therefore $g_n(f(c))=b$ for all but finitely many $n$. Since $C$ is finite, we get $(g_n\circ f)(C) \subseteq B_r$ for all but finitely many $n$.

$(3) \Rightarrow (2)$. Let $c \in C$.  Let $G \preceq G^*$ and $f \in Mon(G^*/A)$.  We claim that $f(c) \in B_r$, so we can  take $B(C)=B_r$. Let $(g_n : f(G) \rightarrow G)_{n \in \mathbb N}$ be a sequence strongly converging to $f(G)$; its existence is assured by Lemma \ref{thm:existence-of-converging-sequence}. So, there exists  $b \in B_r$ such that  $g_{n_k}(f(c))=b$  for some subsequence $(n_k)_{k \in \mathbb N}$. Therefore, by property \ref{item:moves-who-is-outside} of definition \ref{defi:strongly-converges}, we have
$f(c)=b$. Hence $f(C) \subseteq  B_r$ as claimed.

$(2) \Rightarrow (1)$. We suppose that (1) does not hold and we show that (2) does not hold. Let $c\in C \setminus acl^{\exists}(A)$. Then, any existential formula $\phi(x)\in tp^{\exists}(c/A)$ has infinitely many realizations. Define the theory $T(d)=Diag_{el}(G)\cup\{\phi(d),d\not=g_i ; \phi\in tp^{\exists}(c/A),i\in\mathbb N\}$, where $(g_i)_{i\in\mathbb N}$ is an enumeration of the elements of $G$. As $T(d)$ is finitely consistent, there exists an elementary extension $G \preceq G'$ such that $G' \models T(d)$, $d\in G'\setminus G$ and $tp^{\exists}(c/A)\subseteq tp^{\exists}(d/A)$. By Proposition \ref{mono-exten} (2)  there exist an elementary extension $G' \preceq G^*$ and $f\in Mon(G^*/A)$ such that $f(c)=d$.  Hence (2) is not true and this ends the proof. \qed

\bigskip
\noindent \textbf{Proof of Proposition \ref{prop-existent-const}.} Let $D$ be a finite generating set of $K$. Since $acl^\exists(A) < K$ we have $D \not \subseteq acl^\exists(A)$. Hence, using the
equivalence of points 1 and 3  of Lemma \ref{thm:equivalence-acl-ex-extension-converging}, we have:

$(*)$ For any $r \geq  0$ there exist an elementary extension $G^*$ of $G$, a monomorphism  $f \in Mon(G^*/A)$ and a sequence $(g_n: f(G) \rightarrow G)_{n \in \mathbb N}$ strongly converging to $f(G)$, such that $\max_{d \in D}\lh{(g_n\circ f)(d)} \geq r$ for some subsequence $(n_k)_{k \in \mathbb N}$.

Write $K \setminus \{1\}$ as an increasing sequence of finite subsets $(C_i)_{i \in \mathbb N}$.  Enumerate the elements of $acl^\exists(A)$: $acl^\exists(A)=(b_i)_{i \in \mathbb N}$. Let $B_{r(i)}$ be the ball witnessing point 3 of Lemma \ref{thm:equivalence-acl-ex-extension-converging} for $b_i$.

\begin{claim}\label{claim:existence-hom-bounding-acl}
For any $m \in \Bbb N$ there exists a homomorphism $h_m :K \rightarrow G$ satisfying the following
properties:
\begin{enumerate}
 \item\label{item:inject} $1 \not \in h_m(C_m)$;
 \item\label{item:divergence-generators} $\max_{d \in D}\lh{h_m(d)}\geq m$;
 \item\label{item:convergence-acl} $h_m(b_i) \subseteq B_{r(i)}$ for $0 \leq i \leq m$.
\end{enumerate}
\end{claim}

\textit{Proof. }Let $m \in \Bbb N$.  Let $f \in Mon(G^*/A)$ and let $(g_n: f(G) \rightarrow G)_{n \in \mathbb N}$ be the sequence witnessing $(*)$ for $m$.  Since  $(g_n: f(G) \rightarrow G)_{n \in \mathbb N}$ strongly converges to $f(G)$ we have

$$1 \not \in (g_n \circ f)(C_m)$$

for all but finitely many $n$.

Since $b_i \in acl^\exists(A)$, by the equivalence of points 1 and 3 of Lemma \ref{thm:equivalence-acl-ex-extension-converging} we have for any $0 \leq i \leq m$,
$$(g_n\circ f)(b_i) \subseteq B_{r(i)}$$

for all but finitely many $n$.

So, by taking $n_k$ large enough, we obtain:

\begin{enumerate}
 \item $1 \not \in (g_{n_k} \circ f)(C_m)$;
 \item $\max_{d \in D}\lh{(g_{n_k} \circ f)(d)}\geq m$;
 \item $(g_{n_k} \circ f)(b_i) \subseteq B_{r(i)}$ for $0 \leq i \leq m$.
\end{enumerate}

Let $h_m = g_{n_k} \circ f \restr K$. Then $h_m$ is the desired homomorphism and this ends the proof of the Claim. \qed

By point \ref{item:divergence-generators} of the above claim and finiteness of balls of finite radius, we can extract a subsequence $(h_{m_n})_{n \in \mathbb N}$ of pairwise distinct homomorphisms. Thus, we may assume that the initial sequence consists of pairwise distinct homomorphisms.  We are left to show that the sequence $(h_m : K \rightarrow G)_{m \in \mathbb N}$ satisfies the required properties.
By point \ref{item:inject} of Claim \ref{claim:existence-hom-bounding-acl}, the sequence is stable and has a trivial stable kernel. Let $b \in acl^\exists(A)$. Then there exists $p$ such that $b=b_p$. Hence for any $m \geq p$ we have $h_m(b) \in B_{r(p)}$, thus the sequence bounds $acl^\exists(A)$ in the limit.  Therefore, the sequence satisfies all the required properties, so this ends the proof.\qed

\bigskip
To prove Theorem \ref{thm-const-exis} we  need the following result of Takahasi.

\begin{prop} \label{descending-chain-free-groups} \emph{\cite{takahasi}} Let $F$ be a free group of finite rank and let $(L_i |i \in \mathbb N)$ be a descending chain of subgroups with bounded rank. Then $\bigcap_iL_i$ is a free factor of $L_n$ for  all but finitely many $n$.  \qed
\end{prop}

\bigskip
\noindent \textbf{Proof of Theorem \ref{thm-const-exis}.} Define a descending sequence $(L_i|i \in \mathbb N)$ of subgroups of $F$ with bounded rank and containing $acl^\exists(A)$ as follows. Let $L_0=K$.  Suppose that $L_{i}$ is defined.  If $L_{i}=acl^\exists(A)$ then this terminates the sequence; put $L_j=L_i$ for any $j \geq i$.  If $L_i$ is freely $acl^\exists(A)$-decomposable, then set $L_{i+1}$ to be the free factor of $L_i$ containing $acl^\exists(A)$ and which is  freely $acl^\exists(A)$-indecomposable.  So, suppose that $acl^\exists(A)<L_i$ and $L_i$ is freely $acl^\exists(A)$-indecomposable. By Proposition \ref{prop-existent-const} there exists a stable sequence of pairwise distinct homomorphisms $({h_n: L_i \rightarrow F})_{n \in \mathbb N}$ with trivial stable kernel and which bounds $acl^\exists(A)$ in the limit. Hence by Theorem \ref{relative-splitting}  $L_i$ admits a nontrivial  cyclic  splitting  relative to $acl^\exists(A)$.  Then, set $L_{i+1}$ to be the vertex group containing $acl^\exists(A)$.

We claim that the sequence terminates. Suppose for a contradiction that it does not terminate. Then   we get an  infinite  sequence  $(L_i| i \in \mathbb N)$ such that:

$(i)$ $acl^\exists(A)\leq L_i$,

$(ii)$ $rk(L_i) \leq rk(K)$ (properties of free groups, this can be proved   using  \cite{ould-not-cyc} as in Corollary \ref{acl-compress}),

$(iii)$ $L_{i+1} < L_i$.

By Proposition \ref{descending-chain-free-groups}, $\bigcap_{i}L_i$ is a free factor of $L_i$ for all but finitely many $n$. Hence, for all but finitely many $n$, $L_n$ is freely decomposable with respect to $acl^\exists(A)$; a contradiction with the construction of the sequence.  Therefore the sequence terminates, as claimed. Let $L_p$ be the last term in the sequence. Then by construction $acl^\exists(A)=L_p$. We conclude that $K$ is constructible from $acl^\exists(A)$. \qed

As in the case of the algebraic closure, as a consequence we have the following result:

\begin{cor} Let $F$ be a free group of finite rank and $A$ a nonabelian subgroup of $F$. Then $acl^\exists(A)$ is compressed.
\end{cor}

\proof The proof is identical to that of Corollary \ref{acl-compress} by using Theorem \ref{thm-const-exis} instead of Theorem \ref{thme}. \qed

\section{The algebraic closure in the JSJ-decomposition}

In this section we study the link between the algebraic closure and the JSJ-decomposition and we prove Theorem \ref{thm-princip-3}.  We start with the following  lemma.

\begin{lem}\label{lem-acl-factor} Let $G$ be a torsion-free CSA group whose abelian subgroups are cyclic.  Suppose that $G=G_1*G_2$ with $A \leq G_1$. Then $racl_{G}(A) \leq racl_{G_1}(A)$.

\end{lem}

\proof We show first that $racl_{G}(A) \leq G_1$.  We suppose that $g \not \in G_1$ and we find a sequence $(f_n)_{n \in \Bbb N}$ in $Aut(G/A)$ such that the orbit $\{f_n(g); n \in \Bbb N\}$ is infinite; this will prove that $g \not \in racl(A)$. Depending whether $G_2$ is abelian or not, we will treat the two cases separately.  First suppose that $G_2$ is abelian. Then $G_2$ is cyclic; let $t$ be a generating element. Let $\alpha \in G_1$ be nontrivial.  Then, let $(f_n)_{n \in \Bbb N}$ be the sequence of automorphisms of $G$ defined by being the identity on $G_1$ and sending $t$ to $\alpha^n t$. Since $g \not \in G_1$,  $g$ has a normal form $g_0 t^{\epsilon_0} g_1 \cdots g_r t^{\epsilon_r} g_{r+1}$ where $g_i \in G_1$,  $\epsilon_i=\pm 1$ and if $g_i=1$ then $\epsilon_i+\epsilon_{i+1}\neq 0$. If $f_n(g)=f_m(g)$ with $n \neq m$ then a calculation with normal forms shows that $\alpha^{n-m}=1$ which is a contradiction with torsion-freeness of $G$. Hence the orbit $\{f_n(g); n \in \Bbb N\}$ is infinite, as required.

Suppose now that $G_2$ is nonabelian. Since $g \not \in G_1$, $g$ has a normal form $g=g_1 \cdots g_r$,  $r \geq 2$.  Let $g_l \in G_2$ appear in the normal form of $g$. Since $G_2$ is nonabelian and CSA, there exists  an element $\alpha \in G_2$ such that $[g_l, \alpha]\neq 1$. Then, let $(f_n)_{n \in \Bbb N}$ be the sequence of automorphisms of $G$ defined by being  identity on $G_1$ and conjugation by $\alpha^n$ on $G_2$.  If $f_n(g)=f_m(g)$ with $n \neq m$, then a calculation with normal forms shows that $[\alpha^{n-m}, g_l]=1$ which is a contradiction, as $G$ is commutative transitive and $[g_l, \alpha]\neq 1$. Hence the orbit $\{f_n(g); n \in \Bbb N\}$ is infinite, as required.

 Now we show that $racl_G(A) \leq racl_{G_1}(A)$. Let $b \in racl_G(A)$ and suppose that $b \not \in racl_{G_1}(A)$. Then the orbit $\{f(b)| f \in Aut(G_1/A)\}$ is infinite; since each element of $Aut(G_1/A)$ has a natural extension to $G$, the orbit  $\{f(b)| f \in Aut(G/A)\}$ is also infinite, which is a contradiction. \qed

\begin{lem} Let $\G$ be a torsion-free hyperbolic group and $A$ a nonabelian subgroup of $\G$. Suppose that $\G=\G_1*\G_2$ with $A \leq \G_1$ and $\G_1$ is freely $A$-indecomposable. Then $racl_{\G}(A)=racl_{\G_1}(A)$.
\end{lem}

\proof By Lemma \ref{lem-acl-factor}, we have $racl_\G(A) \leq racl_{\G_1}(A)$; thus it remains to show  that $racl_{\G_1}(A) \leq racl_{\G}(A)$.

 Let $f \in Aut(G/A)$. We claim that $f \restr \G_1 \in Aut(\G_1/A)$.  By Grushko-Kurosh theorem, $f(\G)$ has a decomposition $$f(\G_1) =\G_1^{g_1} \cap f(\G_1) *\dots *\G_1^{g_p} \cap f(\G_1)*\G_2^{h_1} \cap f(\G_1) *\dots *\G_2^{h_q} \cap f(\G_1)*F,$$ where $F$ is a free group.  Since $A \leq f(\G_1)$ we have $g_i=1$ for some $i$ and $A \leq \G_1 \cap f(\G_1)$ and this last group is a free factor of $f(\G_1)$. Since $\G_1$ is freely $A$-indecomposable, we conclude that    $\G_1 \cap f(\G_1)=f(\G_1)$, thus $f(\G_1) \leq \G_1$.  If $f \restr \G_1$ is not an automorphism, then by Corollary \ref{cor-mono}, $\G_1$ is freely $A$-decomposable, which is a contradiction. Hence $f \restr \G_1 \in Aut(\G_1/A)$, as claimed.

 Therefore, if the orbit $\{f(b)| f\in Aut(\G_1/A)\}$ is finite then the orbit $\{f(b)| f\in Aut(\G/A)\}$ is finite as well, which proves $racl_{\G_1}(A) \leq racl_{\G}(A)$. \qed

\begin{prop}  \label{prop-vertex}  Let $G$ be a torsion-free CSA group  and $A$ a subgroup of $G$.  Let $\Lambda$ be  an abelian splitting of $G$ relative to $A$ and suppose that each edge group is  maximal abelian in its endpoints vertex groups. If $G(A)$ is the vertex group containing $A$ then $racl(A) \leq G(A)$ and in particular $acl(A) \leq G(A)$.
\end{prop}

\proof  As in the proof fo Lemma \ref{lem-acl-factor}, we are going to show that if $g \not \in G(A)$ then there exists a sequence $(f_n)_{n \in \Bbb N}$ in $Aut(G/A)$ such that the orbit $\{f_n(g); n \in \Bbb N\}$ is infinite; which proves that $g \not \in racl(A)$.  Let $g \not \in G(A)$.

Write $\Lambda=(\mathcal G(V,E), T, \phi)$. To simplify, identify $G$ with $\pi(\mathcal G(V,E), T)$. Enumerate the edges which lie outside $T$ as $e_1, \dots, e_p$. Let $\mathcal G_i(V,E_i)$ be the graph of groups obtained by deleting $e_i$. Hence $G$ is an HNN-extension of the fundamental group $G_i=\pi(\mathcal G_i(V, E_i),T)$.

Suppose that $g \not \in G_i$.  Write $G=\<G_i, t| C^t=\varphi(C)\>$. Let $c \in C$ be nontrivial. In this case  let $(f_n)_{n \in \Bbb N}$  be the sequence of Dehn twists around $c^n$, that is $f_n$ is  defined by being  identity on $G_i$ and sending $t$ to $c^n t$. As in the previous lemma,  $g$ has a normal form $g_0 t^{\epsilon_0} g_1 \dots g_r t^{\epsilon_r} g_{r+1}$; if $f_n(g)=f_m(g)$, with $n \neq m$,   we find $\alpha^{n-m}=1$, a contradiction with torsion-freeness of $G$. This shows that  the orbit $\{f_n(g); n \in \Bbb N\}$ is infinite, as required.

Suppose that $g \in \cap_{1 \leq i \leq p}G_i$. Note that $\cap_{1 \leq i \leq p}G_i$ is the fundamental group $L$ of the graph of groups $\mathcal G(V, E')$ obtained by deleting all the edges $e_1, \dots, e_p$, relative to the maximal subtree $T$. Let $f_1, \dots, f_q$ be the edges incident to $G(A)$.     Hence, for each $1 \leq i \leq q$, $L$ can be written as an amalgamated free product $L=L_{i1}*_{C_i}L_{i2}$ where $L_{i1}$ and $L_{i2}$ are the fundamental groups of the connected components of the graph obtained by deleting $e_i$ and $G(A) \leq L_{i1}$.

Since $g \not \in G(A)$, there exists $1 \leq i \leq q$ such that $g \not \in L_{i1}$. We claim that there exists a sequence $(f_n)_{n \in \Bbb N}$ in $Aut(L/A)$ such that the orbit $\{f_n(g); n \in \Bbb N\}$ is infinite and such that the restriction of each $f_n$ on any edge group of our initial graph of groups $\mathcal G(V,E)$ is a conjugation by an element of $L$.

Define the sequence $(f_n)_{n \in \Bbb N}$ similarly as in the previous case of HNN-extensions and in  Lemma \ref{lem-acl-factor} above.  Since $g \not \in L_1$, $g$ has a normal form $g=g_1 \cdots g_r$,  $r \geq 2$. Let $g_l \in L_{i2}$ appear in the normal form of $g$.    Let $c \in C$ be nontrivial. In this case let $(f_n)_{n \in \Bbb N}$ be the sequence of Dehn twists around $c^n$;  that is $f_n$ is  defined by being  identity on $L_{i1}$ and  conjugation
by $c^n$ on $L_{i2}$.  If $f_n(g)=f_m(g)$ with $n \neq m$, then a calculation with normal forms shows that $[c^{n-m}, g_l]=1$, thus $[g_l, c]=1$.  Since $C_i$ is maximal abelian, we get $g_l \in C_i$; a contradiction. Hence the orbit $\{f_n(g); n \in \Bbb N\}$ is infinite and the restriction of each $f_n$ on each edge group of $\mathcal G(V,E)$ is a conjugation by an element of $L$, as required.

Each $f_n$ has a standard extension $\hat f_n$ to $G$; thus the sequence $(\hat f_n)_{n \in \Bbb N}$ is a sequence from $Aut(G/A)$ with the orbit $\{\hat f_n(g); n \in \Bbb N\}$ infinite, as required. \qed

\begin{prop} \label{thm-place-of-acl-in-splitting}  Let $\G$ be a torsion-free hyperbolic group and let $A$ be a nonabelian subgroup of $\G$.  Then $racl(A)$ coincides with  the vertex group containing $A$ in the generalized malnormal cyclic JSJ-decomposition of $\G$ relative to $A$.
\end{prop}

\proof

Write $\G=\G_1*\G_2$ with $A \leq \G_1$ and $\G_1$ freely $A$-indecomposable. By Lemma \ref{lem-acl-factor}, $racl_{\G}(A) =racl_{\G_1}(A)$; thus we must show that $racl_{\G_1}(A)$ is the vertex group containing $A$ in the cyclic malnormal JSJ-decomposition of $\G_1$ relative to $A$.   Let $G(A)$ be the vertex group containing $A$.

By  Theorem  \ref{finite-index-auto}, there exists a finite number of automorphisms   $f_1, \dots, f_l $ of $\Gamma_1$  such that for any $f \in Aut(\G_1/A)$,  there exists a modular automorphism $\sigma \in Mod(\G_1/A)$ such that $f =f_i \circ \sigma$ for some $i$.

Let $b \in G(A)$. By Lemma \ref{mod-preservation} any automorphism $\sigma \in Mod(\G_1/A)$ fixes the vertex group containing $A$ in the JSJ-decomposition of $\G_1$ relative to $A$. We see that this last property is steal true for the vertex group $G(A)$. Since any $\sigma \in Mod(\G_1/A)$ fixes $G(A)$ pointwise  for any automorphism $f \in Aut(\G_1/A)$ we have $f(b) \in \{f_1(b), \dots, f_l(b)\}$. Thus $b \in acl_{\G_1}(A)$ and $G(A) \leq acl_{\G_1}(A)$.  The inverse inclusion  follows from Proposition \ref{prop-vertex} and properties of the malnormal JSJ-decompositions stated in Theorem \ref{thm-JSJ-hyp}.  \qed

\bigskip
In the case of free groups, we have a  bit more.

\begin{thm} Let $F$ be a free group of finite rank and let $A$ be a nonabelian subgroup of $F$. Then $acl(A)$ coincides with  the vertex group containing $A$ in the generalized malnormal cyclic JSJ-decomposition of $F$ relative to $A$.
\end{thm}

\proof Write $F=F_1*F_2$ with $A \leq F_1$ and $F_1$ freely $A$-indecomposable.  Since $F_1 \preceq F$, $acl_{F_1}(A)=acl_{F}(A)$. Let $G(A)$ be the vertex group containing $A$ in the cyclic malnormal  JSJ-decomposition of $F_1$ relative to $A$. By Proposition \ref{prop-vertex} and properties of JSJ-decompositions stated in Theorem \ref{thm-JSJ-hyp}, we have $acl(A) \leq G(A)$; thus it remains to show that $G(A) \leq acl(A)$. Let $c \in G(A)$ and let $(\bar d_1, \bar d_2)$ be  a tuple generating $F_1$ with $\bar d_1$ generating $G(A)$. Then $c=w(\bar d_1)$ for some word $w$.

By Theorem \ref{thme}  $acl(A)$ is finitely generated; let $\bar b$  be a finite generating set of $acl(A)$.
Let $\varphi(\bar x, \bar y)$ be the formula given by the Proposition \ref{formula-homogeneity} with respect to the generating tuple $(\bar d_1, \bar d_2)$ and to the tuple $\bar b$; that is for any endomorphism $f$ of $F_1$, if $F_1 \models \varphi(f(\bar d_1), f(\bar d_2))$ and $f$ fixes $\bar b$ then $f$ is an automorphism.

By equational noetherianity, there exists a finite system $S(\bar x, \bar y)$ of equations such that for any $(\bar \alpha, \bar \beta)$ if $F_1 \models S(\bar \alpha, \bar \beta)$ then the map which sends $(\bar d_1, \bar d_2)$ to $(\bar \alpha, \bar \beta)$ extends to an homomorphism.

Let $\bar v(\bar x)$ be a tuple of words such that $\bar b=\bar v(\bar d_1)$.

Let
$$
\psi(z, \bar b):= \exists \bar x \exists \bar y (\varphi(\bar x, \bar y)\wedge z=w(\bar x)\wedge S(\bar x, \bar y) \wedge \bar b =\bar v(\bar x)).
$$

We claim that $\psi(z, \bar b)$ has only finitely many realizations in $F_1$. Indeed, if
$$
F_1 \models \psi(c', \bar b):= \exists \bar x \exists \bar y (\varphi(\bar x, \bar y)\wedge c'=w(\bar x)\wedge S(\bar x, \bar y) \wedge \bar b =\bar v(\bar x)),
$$
then there exists an automorphism $f$ fixing $acl(A)$ pointwise and sending $c$ to $c'$. By Proposition \ref{thm-place-of-acl-in-splitting}  $G(A)=racl(A)$,  thus  the set $\{f(c)| f \in Aut(F_1/A)\}$ is finite. Hence  $\psi(z, \bar b)$ has only finitely many realizations as claimed.  Thus $c \in acl(acl(A))=acl(A)$ as required. \qed

\section{The algebraic closure \& the definable closure}

Putting all the pieces together, in this section we are ready  to give the relation between  algebraic closure and  definable closure.

\begin{thm} \label{lem8} Let $F$ be a free group of finite rank and $A$ a nonabelian subgroup of $F$. Then $dcl(A)$ is a free factor of $acl(A)$.  Similarly, $dcl^\exists(A)$ is a free factor of $acl^\exists(A)$.
\end{thm}

We need the following theorem of Dyer and Scott.

\begin{thm} \emph{ \cite[Proposition 5.3, Ch I]{LyndonSchupp77} \cite{dyer-scott} } \label{thm-dyer-scott}Let $F$ be a free group of finite rank and let $f$ be an automorphism of $F$ of finite order. Then the set of elements of $F$ fixed by $f$ is a free factor of $F$. \qed
\end{thm}

\textit{Proof of Theorem \ref{lem8}. }By Theorem \ref{thme}, $acl(A)$ is finitely generated. Hence, by Grushko-Kurosh theorem, $acl(A)$ has a free decomposition $acl(A)=K*L$,  such that $K$ contains $acl(A)$  and  it is  freely $acl(A)$-indecomposable. We claim that $K=dcl(A)$. Suppose for a contradiction that $dcl(A)<K$ and let $a \in K \setminus dcl(A)$.

\bigskip
\noindent \textbf{Claim 1.} \textit{There exists an automorphism $h$ of $acl(A)$, of   finite order and fixing pointwise $dcl(A)$,  such that $h(a)\neq a$.}

\proof

Since $a \in acl(A)\setminus dcl(A)$,  there exists a formula $\psi(x)$, with parameters from $A$,  such that $\psi(F)$ is finite, contains $a$ and is not a singleton. We claim that there exists $b \in acl(A)$ such that $tp(a|A)=tp(b|A)$ and $a \neq b$.   Set $\psi(F)=\{a, b_1, \dots, b_m\}$ and suppose towards a contradiction that $tp(a|A) \neq tp(b_i|A)$ for all $1 \leq i \leq m$. Thus, for every $1 \leq i \leq m$, there exists a formula $\psi_i(x)$, with parameters from $A$, such that $\psi_i \in  tp(b_i|A)$ and $\neg \psi_i \in tp(a|A)$.  Thus the formula $\psi(x) \wedge \neg \psi_1(x) \wedge \dots \wedge \neg \psi_m(x)$ defines $a$; a contradiction.

Hence, let $b \in \psi(F)$  such that $a \neq b$ and $tp(a|A)=tp(b|A)$. By Proposition \ref{mono-exten}, there exist an elementary extension $F^*$ of $F$ and $f \in Aut(F^*/A)$ such that $f(a)=b$.  Let $h$ be the restriction of $f$ to $acl(A)$. We claim that $h$ has the required properties.

Since $h$ is  restriction of $f$, we get $h(acl(A)) \leq acl(A)$. Let $b \in acl(A)$ and let $\psi_b(x)$ be a formula, with parameters from $A$, such that $\psi_b(F)$ is finite and  contains $b$. Then $h(\psi_b(F)) \leq \psi_b(F)$ and since $\psi_b(F)$ is finite and $h$ is injective we get $h(\psi_b(F)) = \psi_b(F)$. Thus $h$ is surjective and in particular $h$ is an automorphism of $acl(A)$. Moreover, since for any $n$, $h^n$ is an automorphism of $acl(A)$ and $h^n(\psi_b(F)) = \psi_b(F)$, there exists $n \in \mathbb N$ such that $h^n$ fixes  $\psi_b(F)$ pointwise.

Let $\{b_1, \dots, b_m\}$ be a finite generating set of $acl(A)$. Hence, we get $n_1, \dots, n_m$ such that $h^{n_i}(b_i)=b_i$. Therefore $h^{n_1\cdots n_m }(x)=x$ for any $x \in acl(A)$, thus $h$ has finite order. This completes the proof of the claim. \qed

\smallskip
Let $h$ be the automorphism given by the above claim. We claim that $h(K)=K$. We have $h(K) \leq acl(A)$ and by Grushko-Kurosh theorem
$$h(K)=h(K) \cap K^{g_1}* \dots *h(K) \cap K^{g_n}*h(K) \cap L^{h_1}* \dots *h(K) \cap L^{h_m}*D,$$
where $D$  is a free group.  Since $dcl(A) \leq K \cap h(K)$,  it follows that $g_i=1$ for some $i$.  Since $K$ is $dcl(A)$-freely indecomposable, we find that $h(K) = h(K) \cap K$, thus $h(K) \leq K$.  In particular $h(a) \in K$.

If $h(K)<K$,  then $K$ is freely $dcl(A)$-decomposable  by Corollary \ref{cor-mono}; a contradiction.   Hence $h(K)=K$.

Since $h$ is a nontrivial automorphism of $K$ of finite order, by Theorem \ref{thm-dyer-scott} $K$ is freely $dcl(A)$-decomposable;  a contradiction. Hence in each case we get a contradiction. Therefore $dcl(A)= K$ as required.

\smallskip
Concerning the existential closure, the proof follows the same method. We only give  a sketch of it by detailing the points where the proof is different.  As above, by Theorem \ref{thm-const-exis} instead of Theorem \ref{thme},   $acl^{\exists}(A)$ is finitely generated; hence we get a free decomposition $acl^{\exists}(A)=K *L$, with $dcl^{\exists}(A) \leq K$ and $K$ is freely $acl^{\exists}(A)$-indecomposable. We let $a \in K \setminus dcl^{\exists}(A)$.  As before, we also have the following.

 \bigskip
\noindent \textbf{Claim 2.}
\textit{ There exists an   automorphism $h$ of $acl^{\exists}(A)$, of   finite order and fixing $dcl^{\exists}(A)$ pointwise,   such that $h(a)\neq a$.}

\proof

The unique point, where the proof here is different, is the use of monomorphisms of an elementary extension rather than automorphisms. Since $a \in acl^{\exists}(A)\setminus dcl^{\exists}(A)$,  there exists an existential  formula $\psi(x)$, with parameters from $A$,  such that $\psi(F)$ is finite, contains $a$ and is not a singleton. The claim here is  that there exists $b \in acl^{\exists}(A)$ such that $tp^{\exists}(a|A) \subseteq tp^{\exists}(b|A)$ and $a \neq b$.  The details are similar and left to the reader. Then, by Proposition \ref{mono-exten}, there exists a monomorphism of an elementary extension $F^*$ of $F$ fixing  $dcl^{\exists}(A)$ pointwise such that $f(a)=b$.  Let $h$  be the restriction of $f$ to $acl^{\exists}(A)$. The rest of the proof works exactly as in Claim 1 and  is left to the reader.
\qed

Also the remaining claims work as in the previous case. This ends the proof of the theorem. \qed

\smallskip
Next theorem is a detailed version of Theorem \ref{thm-princip-5}. We  denote by $End(F/A)$ the set of endomorphisms of $F$ fixing $A$ pointwise.

\begin{thm}\label{thm:counterexample}
Let $A_0$ be a finite set (possibly empty) and
$$A=\<A_0,a,b,u|\>, \;H=A*\<y|\>,$$
$$v=aybyay^{-1}by^{-1}, F=\<H,t|u^t=v\>.$$

Then $F$ is a free group of rank $|{A_0}|+4$ and the following properties hold.

\medskip
$(1)$ If $f \in End(F/A)$  then $f \in Aut(F|A)$, and if $f \restr H \neq id_H$ then $f(y)=y^{-1}$.

$(2)$ $acl(A)=acl^{\exists}(A)=H$.

$(3)$ $dcl(A)=dcl^{\exists}(A)=A$.
\end{thm}

\proof

Clearly  $F$ is a free group of rank $|{A_0}|+4$. We suppose $(1)$ and  we show $(2)$ and $(3)$. Clearly we have
$$
A \leq acl^{\exists}(A) \leq acl(A) \leq racl(A),
$$
and since the subgroups generated by $u$ and $v$ respectively are malnormal in $H$, by Proposition \ref{prop-vertex} we have $racl(A) \leq H$. Thus to show $(2)$ it is sufficient to show that $y \in acl^\exists(A)$. Let
$$
\varphi(z):=\exists \alpha  (u^{\alpha}=azbzaz^{-1}bz^{-1}).
$$
Then $F \models \varphi(y)$. Let $\gamma \in F$ such that $F \models \varphi(\gamma)$. Then the map defined by $f(y)=\gamma$, $f(t)=\alpha$ and  identity on $A$ extends to an endomorphism  of $F$ fixing $A$ pointwise; thus, by $(1)$,  $\gamma= y^{\pm 1}$. Hence $\varphi(z)$ has only finitely many realizations, thus $y \in acl^\exists(A)$ as desired.

We show $(3)$. We  have
$$
A \leq dcl^{\exists}(A) \leq dcl(A) \leq rdcl(A) \leq  racl(A)  \leq H,
$$
thus to show $(3)$ it is sufficient to show that  there exists $g \in Aut(F/A)$ such that for any $ \gamma \in H\setminus A$ we have $g(\gamma) \neq \gamma$.  Let $g$ defined on $H$ by being  identity on $A$ and $g(y)=y^{-1}$.  Then
$$
g(v)=ay^{-1}by^{-1}ayby=ay^{-1}by^{-1}aybyay^{-1}by^{-1}(ay^{-1}by^{-1})^{-1}=d v d^{-1},
$$
where $d= ay^{-1}by^{-1}$. Hence by extending $g$ on $F$ by
$$
g(t)= td^{-1},
$$
we get $g \in Aut(F/A)$ with $g(y)=y^{-1}$.  Now if $\gamma \in H\setminus A$ then $y$ appears in the normal form of $\gamma$, thus $g(\gamma) \neq \gamma$ as required.

\smallskip
The remaining  is devoted to the proof of $(1)$.

\bigskip
\noindent \textbf{Claim 1.} \label{f(y)-in-H} \textit{ Let $f\in End(F/A)$. Then $f(y)\in H$.}
\begin{proof}
Suppose towards a contradiction that $f(y) \not \in H$ and let
\begin{equation}\label{notation-f(y)}
f(y)=\alpha_0t^{\epsilon_0}\ldots\alpha_nt^{\epsilon_n}\alpha_{n+1}
\end{equation}
in normal form where $\alpha_i\in H$ and $\epsilon_i=\pm 1$ for every $i$.
\par
By  definition of $v$ and by HNN relation we have
\begin{equation}\label{HNN-f(v)}
 f(t)^{-1}uf(t)=af(y)bf(y)af(y)^{-1}bf(y)^{-1}.
\end{equation}

Substituting definition (\ref{notation-f(y)})  in equation (\ref{HNN-f(v)}), we have
\begin{equation}\label{fv1}
f(t)^{-1}uf(t)\end{equation}
\begin{equation*}
=a\alpha_0t^{\epsilon_0}\ldots\alpha_nt^{\epsilon_n}\alpha_{n+1}b\alpha_0t^{\epsilon_0}\ldots\alpha_nt^{\epsilon_n}\alpha_{n+1}a\alpha_{n+1}^{-1}t^{-\epsilon_n}\alpha_n^{-1}\ldots
\end{equation*}
\begin{equation*}
 \ldots t^{-\epsilon_0}\alpha_0^{-1}b\alpha_{n+1}^{-1}t^{-\epsilon_n}\alpha_n^{-1}\ldots t^{-\epsilon_0}\alpha_0^{-1}.
\end{equation*}
Compare a cyclically reduced conjugate  for each side of (\ref{fv1}): $u$ for the left side, and a cyclically reduced conjugate $c$ of
\begin{equation}\label{fv1-R}
 \alpha_0^{-1}a\alpha_0t^{\epsilon_0}\ldots\alpha_nt^{\epsilon_n}\alpha_{n+1}b\alpha_0t^{\epsilon_0}\ldots\alpha_nt^{\epsilon_n}\alpha_{n+1}a\alpha_{n+1}^{-1}t^{-\epsilon_n}\alpha_n^{-1}\ldots
\end{equation}
\begin{equation*}
 \ldots t^{-\epsilon_0}\alpha_0^{-1}b\alpha_{n+1}^{-1}t^{-\epsilon_n}\alpha_n^{-1}\ldots t^{-\epsilon_0}
\end{equation*}
for the right side.
\par
There are three subwords in $c$ that could be subject to cancellation.
\begin{enumerate}
 \item One is $\alpha_{n+1}a\alpha_{n+1}^{-1}$.

Note that
 \begin{itemize}
  \item it does not belong to $\q{u}$, since two centralizers of generators cannot be conjugate of each other;
  \item it does not belong to $\q{v}$, since this would imply $\alpha_{n+1}a\alpha_{n+1}^{-1}=v^p$ (as $v$ is root-free), and $v^p$ is  cyclically reduced, while a cyclically conjugate of $\alpha_{n+1}a\alpha_{n+1}^{-1}$ is $a$.
 \end{itemize}
 \item The other two subwords are $\alpha_{n+1}b\alpha_0$ and $\alpha_0^{-1}b\alpha_{n+1}^{-1}$.
\par
If the first one is in $\q{u}$, the second
\begin{itemize}
 \item cannot be in $\q{u}$, because their product $\alpha_{n+1}b^2\alpha_{n+1}^{-1}$ should be in $\q{u}$, but it is not, since it is equal to $b^2$ in the Abelianization $H/[H,H]$.
 \item cannot be in $\q{v}$, because their product $\alpha_{n+1}b^2\alpha_{n+1}^{-1}$ should have the form $u^pv^q$, that in the Abelianization is equal to $u^p(a^2b^2)^q$, but, as said above, it is $b^2$.
\end{itemize}
Symmetrically, if the first one is in $\q{v}$, the second
\begin{itemize}
 \item cannot be in $\q{v}$, because their product $\alpha_{n+1}b^2\alpha_{n+1}^{-1}$ should be in $\q{v}$, but it is not, since it is equal to $b^2$ which is different from $a^2b^2=v$ in the Abelianization.
  \item cannot be in $\q{u}$, because their product $\alpha_{n+1}b^2\alpha_{n+1}^{-1}$ should have the form $v^qu^p$, that in the Abelianization is equal to $u^p(a^2b^2)^q$, but, as said above, it is $b^2$.
\end{itemize}
\end{enumerate}
So, suppose $\alpha_{n+1}b\alpha_0$ is in $\q{u}$ or in $\q{v}$, so that we can reduce between the first and the second
occurrence of $f(y)$.
\par
We have the following two cases:
\begin{enumerate}
 \item the reduction procedure stops somewhere, and we are done, since we have some occurrences of $t$ remaining, at least among the first two occurrences of $f(y)$, getting in this way a contradiction (recall that the HNN length of the cyclically reduced conjugate of the left side of the equation (\ref{fv1}) is 0);
 \item the procedure goes on until every $t$ in the first two occurrences of $f(y)$ is cancelled, and we remain with the word
\begin{equation*}
 \alpha_0^{-1}a\alpha_0 d \alpha_{n+1}a\alpha_{n+1}^{-1}t^{-\epsilon_n}\alpha_n^{-1}\ldots
\end{equation*}
\begin{equation*}
 \ldots t^{-\epsilon_0}\alpha_0^{-1}aba\alpha_{n+1}^{-1}t^{-\epsilon_n}\alpha_n^{-1}\ldots t^{-\epsilon_0}.
\end{equation*}
where $d \in \<u\> \cup \<v\>$. This is cyclically reduced, because $\alpha_0^{-1}a\alpha_0d\alpha_{n+1}a\alpha_{n+1}^{-1}=a^2u^p$ or $a^2(a^2b^2)^q$ in the Abelianization, so the above expression  neither belongs to $\<u\>$ nor to $\<v\>$. Thus, also in this case we get a contradiction, since we cannot cancel the remaining occurrences of $t$.
\end{enumerate}
Symmetrically, if $\alpha_0^{-1}b\alpha_{n+1}^{-1}$ belongs to $\q{u} \cup \q{v}$, then at least the occurrences of $t$ in the first two occurrences of $f(y)$ remain, so we get a contradiction as well.
\par
Thus, we can now say that $\lh{f(y)}_{H\!N\!N}=0$, so Claim 1 is proved.

\end{proof}

\bigskip
\noindent \textbf{Claim 2.} \textit{$f(t)\not\in H$.}
\begin{proof}
Suppose that $f(t)=k\in H$ and let $h=f(y)$. Then, from the equation $f(t)^{-1}uf(t)=f(v)$ we have $k^{-1}uk=ahbhah^{-1}bh^{-1}$, an equation in $H$ that in the Abelianization $H/[H,H]$ becomes $u=a^2b^2$, which is not true, so Claim 2 is proved.
\end{proof}

To prove  next claim, we need the following lemma.

\begin{lem} \label{lem1} \emph{\cite{ould-homogeneity}} Let $G=\<H,t|U^t=V\>$ where $U$ and $V$ are cyclic subgroups of $G$ generated respectively by $u$ and $v$. Suppose that:

$(i)$  $U$ and $V$ are   malnormal in $H$.

$(ii)$ $U^h \cap V =1$ for any $h \in H$.

\smallskip
Let $\alpha, \beta \in H$, $s \in G$  such that $\alpha ^s=\beta$,  $|s| \geq 1$. Then one of the following cases holds:

\smallskip
$(1)$   $\alpha= {u^p}^{\gamma}$,     $\beta= {v^p}^{\delta}$,  $s=\gamma^{-1} t \delta$,  where $p \in \mathbb Z$ and $\gamma, \delta \in H$.

\smallskip
$(2)$   $\alpha= {v^p}^{\gamma}$,     $\beta= {u^p}^{\delta}$,  $s=\gamma^{-1} t^{-1} \delta$,  where $p \in \mathbb Z$ and $\gamma, \delta \in H$.  \qed
\end{lem}

\bigskip
\noindent \textbf{Claim 3.} \textit{There exists $\alpha,\beta  \in A$ such that $f(y)=\alpha y^\varepsilon \beta$ where $\varepsilon =\pm 1$.}

\begin{proof}
Since $f(t) \not \in H$ and $f(v) \in H$, by the above lemma $f(v)$ is conjugate to $v$ in $H$.

First of all, $f(y) \not \in A$. Indeed, if $f(y) \in A$ then $f(v) \in A$ which cannot be $H$-conjugate to $v$.

Let
$$
f(y)=h=h_0y^{\epsilon_0}\ldots h_ny^{\epsilon_n}h_{n+1},
$$
where $\varepsilon=\pm 1$ and $h_i \in A$; moreover, if $h_i=1$ then $y^{\varepsilon_{i-1}}y^{\varepsilon_{i}}\neq 1$.

We obtain that $v$ is a $H$-conjugate of
$$
a(h_0y^{\epsilon_0}\ldots h_ly^{\epsilon_l}h_{l+1})b(h_0y^{\epsilon_0}\ldots h_ly^{\epsilon_l}h_{l+1})a(h_0y^{\epsilon_0}\ldots h_ly^{\epsilon_l}h_{l+1})^{-1}b((h_0y^{\epsilon_0}\ldots h_ly^{\epsilon_l}h_{l+1}))^{-1}.
$$
By a similar argument to Claim 1, we get that the unique possibility is that $n=0$.
\end{proof}

\bigskip
\noindent \textbf{Claim 4.} \textit{$f \in Aut(F/A)$.}

\begin{proof} Immediate from the above lemma and Claim 3.
\end{proof}

\bigskip
\noindent \textbf{Claim 5.} Either $f \restr H =id_H$ or $f(y)=y^{-1}$.

\begin{proof}
By Claim 3 and Lemma \ref{lem1}, we know that $f$ conjugates $v$ in $H$ and $f(y)=\alpha y^{\varepsilon} \beta$, where $\varepsilon=\pm 1$. Therefore, by comparison of cyclically reduced words, the word
$$aybyay^{-1}by^{-1}$$
is a cyclic permutation of the word $$\alpha^{-1}a\alpha y^{\epsilon}\beta b \alpha y^{\epsilon}\beta a\beta ^{-1}y^{-\epsilon}\alpha ^{-1}b\beta ^{-1}y^{-\epsilon}.$$
\par

In both cases $\epsilon=+1$ and $\epsilon=-1$, this yields the equations
\begin{itemize}
 \item $\alpha^{-1}a\alpha=a$
 \item $\beta b \alpha=b$
 \item $\beta a\beta ^{-1}=a$
 \item $\alpha ^{-1}b \beta ^{-1}=b.$
\end{itemize}
From the first and the third equations, $\alpha$ and $\beta$ commute with $a$; so $\alpha=a^p$ and $\beta=a^q$.

From the second equation, we have $p=q=0$.

Therefore, if $\epsilon=+1$, then $f \restr H$ is the identity, while, if $\epsilon=-1$, then $f(y)=y^{-1}$. So this last claim and Theorem \ref{thm:counterexample} are proved.
\end{proof}



\bigskip
\noindent Abderezak OULD HOUCINE, \\
Universit\'e  de Mons, Institut de Mathmatique, B\^atiment Le Pentagone,
Avenue du Champ de Mars 6, B-7000 Mons, Belgique. Universit\'e  de Lyon; Universit\'e Lyon 1; INSA de Lyon, F-69621; Ecole Centrale
de Lyon; CNRS, UMR5208, Institut Camille Jordan, 43 blvd du 11 novembre
1918, F-69622 Villeurbanne-Cedex, France.  \\

\noindent \textit{E-mail}:\textrm{ould@math.univ-lyon1.fr}

\bigskip
\noindent Daniele A.G. Vallino, \\
Dipartimento di Matematica \ap Giuseppe Peano\ch, Universit\`a di Torino\\
Via Carlo Alberto, 8, 10121 Torino, Italy.
\noindent \textit{E-mail}:\textrm{vallino@math.univ-lyon1.fr}


\begin{thebibliography}{MVW07}

\bibitem[Bau67]{Baum}
B.~Baumslag.
\newblock Residually free groups.
\newblock {\em {P}roc. {L}ondon {M}ath. {S}oc}, 3(17):402--418, 1967.

\bibitem[BF95]{bestvina-feighn}
M.~Bestvina and M.~Feighn.
\newblock Stable actions of groups on real trees.
\newblock {\em Invent. Math.}, 121(2):287--321, 1995.

\bibitem[BMR99]{Baum-Rem-algebraic}
G.~Baumslag, A.~Myasnikov, and V.~Remeslennikov.
\newblock Algebraic geometry over groups. {I}. {A}lgebraic sets and ideal
  theory.
\newblock {\em J. Algebra}, 219(1):16--79, 1999.

\bibitem[CF04]{Casanovas-weak}
E.~Casanovas and R.~Farr{\'e}.
\newblock Weak forms of elimination of imaginaries.
\newblock {\em MLQ Math. Log. Q.}, 50(2):126--140, 2004.

\bibitem[CG05]{gui-limit} C. Champetier and  V. Guirardel, Limit groups as 
limits of free groups: compactifying the set of free groups, Israel 
Journal of Mathematics, 146 (2005).

\bibitem[CK73]{Chang-Keisler}
C.~C. Chang and H.~J. Keisler.
\newblock {\em Model Theory}.
\newblock North-Holland, Amsterdam, 1973.

\bibitem[DS75]{dyer-scott}
J.~L. Dyer and G.P. Scott.
\newblock Periodic automorphisms of free groups.
\newblock {\em Comm. Algebra}, 3:195--201, 1975.

\bibitem[GL08]{Guirardel-Levitt-tree-cyl}
V.~{Guirardel} and G.~{Levitt}.
\newblock {Trees of cylinders and canonical splittings}.
\newblock {\em ArXiv e-prints}, November 2008.

\bibitem[GL09]{Guirardel-Levitt-JSJ1}
V.~{Guirardel} and G.~{Levitt}.
\newblock {JSJ decompositions: definitions, existence, uniqueness. I: The JSJ
  deformation space}.
\newblock {\em ArXiv e-prints}, November 2009.

\bibitem[GL10]{Guirardel-Levitt-JSJ2}
V.~{Guirardel} and G.~{Levitt}.
\newblock {JSJ decompositions: definitions, existence, uniqueness. II.
  Compatibility and acylindricity}.
\newblock {\em ArXiv e-prints}, February 2010.

\bibitem[Gui08]{Guirardel-action}
V.~Guirardel.
\newblock Actions of finitely generated groups on {$\Bbb R$}-trees.
\newblock {\em Ann. Inst. Fourier (Grenoble)}, 58(1):159--211, 2008.

\bibitem[GW07]{groves-2007}
D.~Groves and H.~Wilton.
\newblock Conjugacy classes of solutions to equations and inequations over
  hyperbolic groups, 2007.

\bibitem[Hod93]{Hodges(book)93}
W.~Hodges.
\newblock {\em Model theory}, volume~42 of {\em Encyclopedia of Mathematics and
  its Applications}.
\newblock Cambridge University Press, Cambridge, 1993.

\bibitem[JS79]{Jaco-Shalen}
H.~Jaco and P.B. Shalen.
\newblock Seifert ﬁbered spaces in 3-manifolds.
\newblock {\em Mem. Amer. Math. Soc}, 21(220), 1979.

\bibitem[KM98a]{Kharla-Mias-Irred1}
O.~Kharlampovich and A.~Myasnikov.
\newblock Irreducible affine varieties over a free group. {I}. {I}rreducibility
  of quadratic equations and {N}ullstellensatz.
\newblock {\em J. Algebra}, 200(2):472--516, 1998.

\bibitem[KM98b]{Kharla-Mias-Irred}
O.~Kharlampovich and A.~Myasnikov.
\newblock Irreducible affine varieties over a free group. {II}. {S}ystems in
  triangular quasi-quadratic form and description of residually free groups.
\newblock {\em J. Algebra}, 200(2):517--570, 1998.

\bibitem[Kro90]{kropholler}
P.H. Kropholler.
\newblock An analogue of the torus decomposition theorem for certain poincar´e
  duality groups.
\newblock {\em Proc. London Math. Soc}, 60(3):503--529, 1990.

\bibitem[KW99]{kapovich-two}
I.~Kapovich and R.~Weidmann.
\newblock On the structure of two-generated hyperbolic groups.
\newblock {\em Math. Z.}, 231(4):783--801, 1999.

\bibitem[LS77]{LyndonSchupp77}
R.~C. Lyndon and P.~E. Schupp.
\newblock {\em Combinatorial group theory}.
\newblock Springer-Verlag, Berlin, 1977.
\newblock Ergebnisse der Mathematik und ihrer Grenzgebiete, Band 89.

\bibitem[Mar02]{Marker}
D.~Marker.
\newblock {\em Model Theory : An Introduction}.
\newblock Springer-Verlag, New york, 2002.
\newblock Graduate Texts in Mathematics.

\bibitem[MTB10]{medvedev-takloo}
A.~Medvedev and R.~Takloo-Bighash.
\newblock An invitation to model-theoretic galois theory.
\newblock {\em ArXiv}, 2010.

\bibitem[MV04]{ventura-martino-compress}
A.~Martino and E.~Ventura.
\newblock Fixed subgroups are compressed in free groups.
\newblock {\em Comm. Algebra}, 32(10):3921--3935, 2004.

\bibitem[MVW07]{miasnikov-ventura-weil-algebraic}
A.~Miasnikov, E.~Ventura, and P.~Weil.
\newblock Algebraic extensions in free groups.
\newblock In {\em Geometric group theory}, Trends Math., pages 225--253.
  Birkh\"auser, Basel, 2007.

\bibitem[OH07]{Ould-equa}
A.~Ould~{H}oucine.
\newblock Limit {G}roups of {E}quationally {N}oetherian {G}roups.
\newblock In {\em Geometric group theory}, Trends Math., pages 103--119.
  Birkh\"auser, Basel, 2007.

\bibitem[OH10]{ould-not-cyc}
A.~Ould~{H}oucine.
\newblock Note on free conjugacy pinched one-relator groups.
\newblock {\em Submited}, 2010.

\bibitem[OH11]{ould-homogeneity}
A.~Ould~{H}oucine.
\newblock Homogeneity and prime models in torsion-free hyperbolic groups.
\newblock {\em Confluentes Mathematici}, 3(1):121--155, 2011.

\bibitem[Per08]{chloe-these}
C.~Perin.
\newblock Plongements \'el\'ementaires dans un groupe hyperbolique sans
  torsion.
\newblock {\em Th\`ese de doctorat, Universit\'e de Caen/Basse-Normandie},
  2008.

\bibitem[Poi83]{poizat-imaginaire}
B.~Poizat.
\newblock Une th\'eorie de {G}alois imaginaire.
\newblock {\em J. Symbolic Logic}, 48(4):1151--1170 (1984), 1983.

\bibitem[PS10]{perin-homo}
C.~Perin and R.~Sklinos.
\newblock Homogeneity in the free group.
\newblock {\em Preprint}, 2010.

\bibitem[RS94]{sela-rips-rigid}
E.~Rips and Z.~Sela.
\newblock Structure and rigidity in hyperbolic groups. {I}.
\newblock {\em Geom. Funct. Anal.}, 4(3):337--371, 1994.

\bibitem[RS97]{sela-JSJ}
E.~Rips and Z.~Sela.
\newblock Cyclic splittings of finitely presented groups and the canonical
  {JSJ} decomposition.
\newblock {\em Ann. of Math. (2)}, 146(1):53--109, 1997.

\bibitem[RW10]{weidmann-equa}
C.~Reinfeldt and R.~Weidmann.
\newblock {M}akanin-razborov diagrams for hyperbolic groups.
\newblock {\em Preprint}, 2010.

\bibitem[Sel97a]{sela-acylin}
Z.~Sela.
\newblock Acylindrical accessibility for groups.
\newblock {\em Invent. Math.}, 129(3):527--565, 1997.

\bibitem[Sel97b]{Sela-hopf}
Z.~Sela.
\newblock Structure and rigidity in ({G}romov) hyperbolic groups and discrete
  groups in rank {$1$} {L}ie groups. {II}.
\newblock {\em Geom. Funct. Anal.}, 7(3):561--593, 1997.

\bibitem[Sel01]{Sela-Diophan1}
Z.~Sela.
\newblock Diophantine geometry over groups. {I}. {M}akanin-{R}azborov diagrams.
\newblock {\em Publ. Math. Inst. Hautes \'Etudes Sci.}, (93):31--105, 2001.

\bibitem[Sel06]{Sela-stab}
Z.~Sela.
\newblock Diophantine geometry over groups viii:stability.
\newblock {\em Preprint}, Available at:http://www.ma.huji.ac.il/~zlil/, 2006.

\bibitem[Sel09a]{sela-imaginaries}
Z.~Sela.
\newblock Diophantine geometry over groups ix: Envelopes and imaginaries.
\newblock {\em ArXiv}, 2009.

\bibitem[Sel09b]{Sela-hyp}
Z.~Sela.
\newblock Diophantine geometry over groups. {VII}. {T}he elementary theory of a
  hyperbolic group.
\newblock {\em Proc. Lond. Math. Soc. (3)}, 99(1):217--273, 2009.

\bibitem[Tak51]{takahasi}
M.~Takahasi.
\newblock Note on chain conditions in free groups.
\newblock {\em Osaka Math.J.}, 3(2):221--225, 1951.

\bibitem[Val11]{vallino-thesis}
D.~Vallino.
\newblock {\em Algebraic and definable closure in free groups}.
\newblock PhD thesis. In preparation, 2011.

\bibitem[Wil06]{wilton-thesis}
H.~Wilton.
\newblock {\em Subgroup separability of limit groups}.
\newblock PhD thesis, Univ. London, 2006.

\end{thebibliography}
\end{document}